\documentclass[final]{aupl}

\usepackage{
amsfonts,
latexsym,
amssymb,
enumerate,
mathrsfs,
graphicx,
centernot,
color,
}

\newcommand{\labbel}[1]{\label{#1} [[{\bf #1}]]}  
\renewcommand{\labbel}{\label}

\newtheorem{theorem}{Theorem}[section]

\newtheorem{prop}[theorem]{Proposition} 
\newtheorem{proposition}[theorem]{Proposition} 
 
\newtheorem{corollary}[theorem]{Corollary}

\newtheorem*{claim*}{Claim}

\newtheorem*{theorem*}{Theorem}
\newtheorem*{proposition*}{Proposition}
\newtheorem*{corollary*}{Corollary}
\newtheorem*{lemma*}{Lemma}
\newtheorem*{scholion*}{Scholion}

\theoremstyle{definition}
\newtheorem{definition}[theorem]{Definition}

\newtheorem{problem}[theorem]{Problem}

\theoremstyle{remark}
\newtheorem{remark}[theorem]{Remark}
\newtheorem{remarks}[theorem]{Remarks}
\newtheorem*{remark*}{Remark}
\newtheorem*{remarks*}{Remarks}
 
\newtheorem{example}[theorem]{Example}

\newtheorem{observation}[theorem]{Observation}

\newtheorem*{observation*}{Observation}

 \allowdisplaybreaks[1]

\numberwithin{equation}{section}

\begin{document}

\title{The strong amalgamation property into union}

\author{Paolo Lipparini} 
\address{Dipartimento di Amalgamatematica\\Viale della  Ricerca
 Scientifica\\Universit\`a di Roma ``Tor Vergata'' 
\\I-00133 ROME ITALY}

\email{lipparin@axp.mat.uniroma2.it}



\date{\today}

\thanks{
Work performed under the auspices of G.N.S.A.G.A. Work 
partially supported by PRIN 2012 ``Logica, Modelli e Insiemi''.
The author acknowledges the MIUR Department Project awarded to the
Department of Mathematics, University of Rome Tor Vergata, CUP
E83C18000100006.}

\begin{abstract}
We consider the situation in which 
some class of structures has the Strong Amalgamation Property (SAP)
with the further requirement that the amalgamating structure can be taken
over the set theoretical union of (the images of) the domains of the
 structures to be amalgamated.
We call this property SAPU. 

The main advantage of SAPU over SAP is that there 
are many preservation theorems showing that we can 
merge different theories with SAPU still obtaining a theory 
with SAPU, hence with SAP.
In particular, we get  SAPU
for various theories with  many
binary relations, each  relation
satisfying any set  of properties chosen among 
transitivity,
 reflexivity,
 symmetry,
antireflexivity,
antisymmetry.
We may also
add unary operations, possibly satisfying
some coarseness, isotonicity and closure conditions.

SAPU is not limited to relational theories:
 the varieties defining the most usual Maltsev 
conditions in universal algebra have SAPU.
Other examples include bounded directoids, order algebras
and various generalizations.
 \end{abstract}

\maketitle

\section{Introduction.} \labbel{intro} 

The Amalgamation Property 
(AP) has found many important applications
in algebra, logic, category theory and, recently,
computer science.
See, e.~g., \cite{BGR,H,GM,GG,J,KMPT,Mac,Ma,MMT}. 

We study theories with the amalgamation property 
and with the further requirement that the amalgamating structure can be taken
over the set-\hspace{0 pt}theoretical union of the (images of the) structures to be amalgamated.
This applies to any universal theory in a purely relational 
language and with AP, in particular,
this is the case for partially ordered sets
(henceforth
\emph{posets}, for short), more generally,
for structures with a binary relation satisfying any 
(fixed in advance) number of the following properties:
transitivity,
 reflexivity,
 symmetry,
antireflexivity,
antisymmetry. See Proposition \ref{pu}. 

The advantage of AP over Union (APU)
is that we can frequently merge different theories with APU,
still obtaining a theory with APU, hence with AP.
For example, this holds for theories with a pair of binary 
relations as above, possibly with the condition asserting
that one relation is coarser than (i. e., contains) another.
See Theorem \ref{due}.
In general, the result fails when three relations are taken into account.

On the other hand, we do have APU when all the relations
under consideration are assumed to be transitive; this applies
to an arbitrary number of relations.
In particular, APU holds for any class of multiposets
with any prescribed set of coarseness relations.
With a few exceptions, APU is  maintained when  unary operations 
preserving one or more relations are added;
in particular, all multiposets with 
operators have APU. See Theorems \ref{mul} and  \ref{gen}. 
Many applications of APU are relatively simple, however, most
of the results are really fine-tuned, in that, just weakening some assumption,
counterexamples can be found, e.~g., Theorem \ref{lopp}(7), Examples \ref{exbin},
\ref{inf}(d), \ref{necess}, \ref{uneed}, \ref{1}, \ref{injbij}, \ref{exch} and   
  Propositions \ref{div+}, \ref{dense}, \ref{bin}, \ref{rt}.

There are results special to APU which generally do not hold for
AP. If $T$ is a theory with APU and 
we add to $T$ universal-existential sentences in which only one variable is bounded
by the universal quantifier, then the resulting theory has still APU.
In particular, if some class of 
partially ordered sets  with a unary operation has APU,
we still have APU if we ask that the operation is
an involution, or that the  operation,
if  order preserving, is a closure operation.
In general, there are plenty of conditions $H$
such that if $\mathcal K$ is a class with APU,
then the subclass of those structures in $\mathcal K$ satisfying
$H$ maintains APU.  
See Section \ref{pres}.
Counterexamples are provided for theories with 
AP  not into Union.

In another direction, studying APU is useful for discovering results
which hold in general for AP.
Dealing with the Strong APU (see below for the definition),
it is almost immediate to show that the union of theories in disjoint
languages with SAPU has still SAPU.
This fact turns out to be true for SAP, as well,
but with a not entirely trivial argument.
Counterexamples exist showing that it is necessary
to deal with the strong variants.

At first sight, the reader might expect that APU
is a phenomenon almost exclusively limited
to relational languages. This is not the case.
On one hand, we can  consider unary
functions, usually getting APU ``almost for free''.
On the other hand, there is a bunch of examples 
of theories with APU in  languages with $n$-ary functions,
 for $n \geq 2$. See Section \ref{binbin} but also
Propositions \ref{T} and \ref{un}.
Not only APU 
has many interesting and useful consequences,
but it is applicable to a number of nontrivial examples.

\section{The Strong Amalgamation Property into
 Union.} \labbel{sec} 

We work 
with classes of   structures 
with finitary
relations and functions. \emph{Structure} and \emph{model}
are synonymous.   
 As usual in model 
theory, equality is considered as a logical symbol,
namely, it can be interpreted in every structure, and it is 
actually interpreted as identity. Under a frequent terminology,
this means that we work
with \emph{normal} models.
In particular, we do not include equality in the symbols
belonging to some language $\mathscr L$, so that when $\mathscr L= \emptyset $,
then $\mathscr L$ is the pure language of identity. 

We do not take explicit position 
on the admissibility or not of structures with empty domain.
Generally, our results hold in both settings;
otherwise, we shall mention the assumptions explicitly. 

An \emph{embedding}
 $\iota$ from some structure $\mathbf A$
into a structure
 $\mathbf B$  
for the same language is an injective function
$\iota:A \to B$ 
such that, for every $a_1,a_2, \dots \in A$,
the following hold: 
\begin{equation*}  
 \iota (f_{\mathbf A}(a_1,a_2, \dots)) =
 f_{\mathbf B} (\iota(a_1), \iota(a_2), \dots),
\end{equation*}   
for every function symbol $f$ in the language, and
\begin{equation}\labbel{R}  
 R_{\mathbf A}(a_1,a_2, \dots)  \quad \text{ if and only if} \quad 
 R_{\mathbf B} (\iota(a_1), \iota(a_2), \dots),
 \end{equation}    
for every relation symbol  $R$ in the language.
 
If we drop the requirement of injectivity 
and weaken   condition \eqref{R}  to
\begin{equation} \labbel{RR}    
 R_{\mathbf A}(a_1,a_2, \dots) \quad \text{implies} \quad
 R_{\mathbf B} (\iota(a_1), \iota(a_2), \dots),
\end{equation}     
we get the weaker notion of a \emph{homomorphism}.
Here we shall consider amalgamation properties with respect to 
embeddings. Were we considering injective homomorphisms,
instead, we would get completely different results.
See Remark \ref{mor} below. 

We shall possibly deal also with constants (= selected elements).
Embeddings and homomorphisms are assumed
to satisfy 
$\iota (c_{\mathbf A})= c_{\mathbf B}$,
for every constant symbol $c$. Subscripts will be dropped when
no risk of confusion might arise.
Full formal details about the above notions can be found
in any textbook on model theory, e.~g., 
\cite{H}. 

If $\mathbf A$ and $\mathbf B$
are structures for the same language
and $A \subseteq B$ (as sets), then
we say that   \emph{$\mathbf A$ is a substructure of $\mathbf B$} 
if the inclusion map from $A$ to $B$ is an embedding
of $\mathbf A$ into $\mathbf B$.
In the above situation we shall write
$\mathbf A \subseteq \mathbf B$.

\begin{definition} \labbel{sapu}
(a) A class $\mathcal K$ of structures  
for the same language has the \emph{amalgamation property}
(AP) if, whenever 
$\mathbf A, \mathbf B, \mathbf C \in \mathcal K$,
 $ \iota _{\mathbf C, \mathbf A} \colon \mathbf C \rightarrowtail \mathbf A$
and 
 $ \iota _{\mathbf C, \mathbf B} \colon \mathbf C \rightarrowtail \mathbf B$
are embeddings, then there is a structure
$\mathbf D \in \mathcal K$ and  embeddings
$ \iota _{\mathbf A, \mathbf D} \colon \mathbf A \rightarrowtail \mathbf D$
and 
 $ \iota _{\mathbf B, \mathbf D} \colon \mathbf B \rightarrowtail \mathbf D$
such that 
$\iota _{\mathbf C, \mathbf A} \circ \iota _{\mathbf A, \mathbf D}
=
\iota _{\mathbf C, \mathbf B} \circ \iota _{\mathbf B, \mathbf D}$. 
Namely, the following diagram can be commutatively completed as requested.
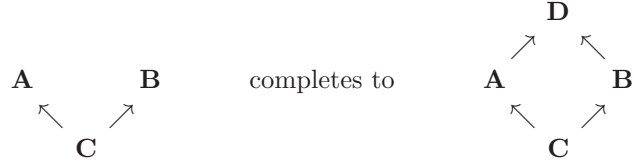
\begin{equation*}
\begin{matrix} 
 \cr
 \cr
\mathbf A \quad \quad \quad \quad \mathbf B \cr
\nwarrow \ \quad\  \nearrow \cr
 \mathbf C
 \end{matrix}   
\qquad \quad  \text{ completes to }  \quad \qquad 
\begin{matrix} 
 \mathbf D \cr
\nearrow \ \quad\  \nwarrow \cr
\mathbf A \quad \quad \quad \quad \mathbf B \cr
\nwarrow \ \quad\  \nearrow \cr
 \mathbf C
 \end{matrix}   
\end{equation*}    

(b) A class of structures  has the \emph{strong amalgamation property}
(SAP) if, under the assumptions in (a), the conclusion can be
strengthened to the effect that the intersection
of the images of $\iota _{\mathbf A, \mathbf D}$
and $\iota _{\mathbf B, \mathbf D}$ 
is equal to the image of 
$\iota _{\mathbf C, \mathbf A} \circ \iota _{\mathbf A, \mathbf D}$
(hence also of $\iota _{\mathbf C, \mathbf B} \circ \iota _{\mathbf B, \mathbf D}$).

(c) A class of structures  has the \emph{(strong) amalgamation property
into union}
(SAPU)  APU if, in addition, 
$\mathbf D$ can be chosen in such a way that
its domain $D$ is the union of the images 
of $\iota _{\mathbf A, \mathbf D}$
and $\iota _{\mathbf B, \mathbf D}$.

Here ``union'' is meant in the set theoretical sense, not
in the model theoretical sense: we are considering the
union of the domains, not of the structures. See below for more details. 

If $T$ is a first-order theory,
we say that $T$ has AP, SAP, APU, SAPU
if the class of models of $T$ has the property.

Formally, we assume that an empty class
of structures shares all the above properties.
This is consistent with standard conventions about
universal quantification over empty domains.
\end{definition}    

\begin{remark} \labbel{iso}    
In our present context, if $\mathcal K$ is closed under isomorphism,
then the above definitions (b) and (c)
can be simplified.

(d) A class $\mathcal K$  closed under isomorphism
has SAP (resp., SAPU) if and only if, whenever 
$\mathbf A, \mathbf B, \mathbf C \in \mathcal K$,
 $ \mathbf C \subseteq  \mathbf A$,
 $ \mathbf C \subseteq  \mathbf B$
and $C=A \cap B$,
then there is a structure
$\mathbf D \in \mathcal K$ such that 
$ \mathbf A \subseteq  \mathbf D$,
$ \mathbf B \subseteq  \mathbf D$
(and, resp., $A \cup B =D$).

We shall sometimes informally refer to a triple
$\mathbf A, \mathbf B, \mathbf C $
as above as a \emph{triple to be amalgamated}, 
for short, a \emph{TBA triple}.

For simplicity, we shall generally work in the simplified setting
described in the previous paragraph, namely, we shall deal with
inclusions as above, rather than with arbitrary embeddings
as in Definition \ref{sapu}.
In particular, we shall always assume that classes of structures are closed under isomorphism.
In any case, the setting in which we  work
 shall always be clear from the context. 
\end{remark}

Results about  AP and SAP appear scattered in the literature,
sometimes in different settings or terminology.
 A survey of results about AP
and related properties appears in 
\cite{KMPT}, where the notions
are also inserted in a general 
categorical framework.
A survey of various applications
of AP to model theory can be found in \cite{H}.

For relational languages, the special case of SAPU
when $\mathbf  D$ can be taken as the model-theoretical
union of $\mathbf A$ and $\mathbf  B$ 
has been considered by various authors,
generally under the name \emph{free amalgamation}.
This is the particular case when
 $R_{\mathbf D} = R_{\mathbf A} \cup R_{\mathbf B}$, 
for every relation symbol $R$ in the language.
For example, see  \cite{B,F2,Mac} and further
references there. 

To the best of our knowledge, the explicit definitions
of APU and  
 SAPU in the general case 
are new when considered  for a whole class of structures.
The first implicit appearance of SAPU
possibly occurs in 
Fra\"\i ss\'e argument \cite[Section 9]{F} 
showing that the class of 
linear orders and the class of
 posets 
have SAP.
Compare also J{\'o}nsson \cite[Lemma 2.3]{J}. 
In another direction, particular situations in which the amalgamating
structure can be taken over $A \cup B$
have been considered in lattice theory.
See \cite[IV, Section 2.3 and VI, Exercise  4.11]{G} and further references there.

The main interest of SAPU comes from the fact that 
there are various methods to join or modify
some theories with SAPU 
 in order to obtain
other theories with  SAPU. 
See  Section \ref{pres} below.
In particular, if we 
merge distinct theories in disjoint languages and having SAPU,
we still obtain a theory with SAPU. See Proposition \ref{T} below.
The result is true also for SAP (Proposition \ref{senzau}),
however, in the case of SAPU the proof is simpler and,
as the main point in the present note,
in certain cases the argument applies also
to theories which are not in disjoint languages. 
See Theorems \ref{mul} and \ref{gen}.  

Notice that there are theories with both SAP and APU
but without SAPU. See Proposition \ref{dense} 
and Example \ref{1a}. Compare also Example \ref{***}.

If $R$ is a binary relation on some set, 
we frequently
write $ a \mathrel {R } b $
in place of $R(a,b)$ or $(a,b) \in R$.   
Moreover, $ a \mathrel {R } b \mathrel { S } c  $
is a shorthand for $ a \mathrel {R } b $ and $b  \mathrel { S } c  $.

In most cases, we shall prove a property somewhat stronger than 
SAPU.

\begin{definition} \labbel{super}
If $\mathcal K$ is a class of structures
closed under isomorphism and  in a language with a 
binary relation symbol $R$,
we say that $\mathcal K$ has  \emph{superSAPU}
(resp., \emph{superSAP})
\emph{with respect to $R$}\/
if, whenever 
$\mathbf A, \mathbf B, \mathbf C \in \mathcal K$
is a triple to be amalgamated as in Remark \ref{iso}, then 
there exists an amalgamating structure $\mathbf D
\in \mathcal K$  witnessing SAPU (resp.,  SAP)  and such that,
for every  $a \in A \setminus B$ and $b \in B \setminus A$, 
  \begin{enumerate}[(i)]
    \item 
if $a \mathrel { R _{\mathbf D}} b $  
then there is $c \in C$ such that  $a \mathrel { R _{\mathbf A}} c 
\mathrel { R _{\mathbf B}} b $, and
\item
if $b \mathrel { R _{\mathbf D}} a $  
then there is $c \in C$ such that  $b \mathrel { R _{\mathbf B}} c 
\mathrel { R _{\mathbf A}} a $.
  \end{enumerate} 

Whenever we speak of superSAP(U), we always assume that at least
one relation $R$ is specified, and superSAP(U) is meant with respect to 
the specified relation(s).
We shall omit the reference to 
the relations, when they are understood.
 \end{definition}    

Usually, the superamalgamation property
(not necessarily into union)  is considered
with respect to  some
ordering relation, see \cite{GM}, but
many applications to algebraic logic
are known even in the case of an arbitrary binary relation,
see \cite{Ma}.
The assumption that
$a \in A \setminus B$ and $b \in B \setminus A$
in the hypothesis of superamalgamation properties
is frequently weakened to 
$a \in A $ and $b \in B $.
We shall use the modified version as in
 Definition \ref{super} in order to simplify statements. 
Of course, the two definitions are equivalent if
$R$ is a reflexive relation (since if, say,
$a \in A \cap B$ and $a \mathrel { R} b$, then we can take $c=a$ in (i),
and similarly  for (ii)).

Another  property related to AP has proven very important in model theory.

\begin{definition} \labbel{jepdef}
The joint embedding property is the special instance of AP
when the structure $\mathbf  C$ is empty.
In detail, a class $\mathcal K$ of structures  
for the same language has the \emph{joint embedding  property}
(JEP) if, for every  
$\mathbf A, \mathbf B  \in \mathcal K$,
 there are a structure
$\mathbf D \in \mathcal K$ and  embeddings
$ \iota: \mathbf A \rightarrowtail \mathbf D$
and 
 $ \kappa: \mathbf B \rightarrowtail \mathbf D$. 
\end{definition}   

In Section \ref{jep} 
we shall present and study variations on JEP in the same spirit
of Definition \ref{sapu}.

\begin{remark} \labbel{jepeq} 
   As mentioned, when we allow the empty structure, JEP is a special instance
of AP. However, there are situations in which it is not convenient to consider
the empty structure or, plainly, such a structure does not exist, e.~g., 
when the language has some constant.

Nevertheless, even assuming AP only for nonempty structures,
we can always reduce ourselves to a situation in which JEP holds.
As well-known, assuming AP  for nonempty structures,
if we set $\mathbf A \sim \mathbf B$
when $\mathbf A $ and $  \mathbf B $
can be joint embedded into some $\mathbf  D$ as in Definition \ref{jepdef},
then $ \sim$ turns out to be an equivalence relation, hence
JEP holds when restricted to each equivalence class.
Clearly, AP is maintained relative to each equivalence class,
so that we have both AP and JEP on each equivalence class.
\end{remark}

If  $\mathcal K$  is a class of finitely generated 
structures in a countable language, a  \emph{Fra\"\i ss\'e limit} 
of $\mathcal K$ is a 
countable ultrahomogeneous structure of age  $\mathcal K$.
 See \cite[Section 7.1]{H} for  details.
A countable (up to isomorphism) hereditary  class $\mathcal K$  with AP and JEP has
a Fra\"\i ss\'e limit \cite[Theorem 7.1.2]{H}.

Throughout, we suppose that all classes of models under consideration
are closed under taking isomorphism. Compare Remark \ref{iso}.

\section{Orders,  binary relations, adding operators.} \labbel{ex}

For relational languages 
the next observation is folklore. 
Here and in similar situations below we point
out the lesser known fact that the argument carries over when
considering also unary function symbols.

\begin{observation} \labbel{obs}
(a) A universal theory in a language without function symbols
of arity $\geq 2$  has APU (resp., SAPU, superSAPU)
if and only if it has AP (resp., SAP, superSAP). 

(b)  More generally, suppose that $\mathcal K$ is a class of structures for
a language without function symbols
of arity $\geq 2$ and suppose further that $\mathcal K$ 
is preserved under taking substructures. Then $\mathcal K$   has
APU (resp., SAPU,  superSAPU)
if and only if $\mathcal K$  has AP (resp., SAP, superSAP). 

(c)  The class of all models for some given language---in 
other words,  any theory without nonlogical axioms---has SAPU and, in the presence of a binary relation,
superSAPU.
\end{observation}   

\begin{proof}
(a) To prove the non trivial implication, 
under the assumptions of (super)(S)APU, take some $\mathbf D$ witnessing
(super)(S)AP  and consider the union $D_1$ of the images of  
$\iota _{\mathbf A, \mathbf D}$
and of $\iota _{\mathbf B, \mathbf D}$.
Since $T$ is universal, then the restriction $\mathbf D_1$
of $\mathbf D$ to $D_1$  
is a model of $T$, thus $\mathbf D_1$ 
witnesses (super)(S)APU. Notice that $\mathbf D_1$
is actually a structure, since it is the union of two structures 
and then unary functions do not send elements of $D_1$
outside.  As far as superSAP(U) is concerned, notice that
the defining property for superSAP(U) speaks only of elements
in the images of $\iota _{\mathbf A, \mathbf D}$
and  $\iota _{\mathbf B, \mathbf D}$, of elements of $A \cup B$,
in the simplified setting of Remark \ref{iso}.  

(b) is proved in a similar way.

(c) The argument
generalizes \cite[9.1]{F}. Relations cause no trouble
and operations can be extended in an arbitrary way in the union,
since
no axiom is prescribed.
\end{proof}

Most results in the present section are taken
from \cite{lop,apuduerel}, where they are stated without the 
specification ``into Union''.
However, the ``U''  follows 
from Observation \ref{obs}
(and, in any case, directly from the original proofs).

See \cite{apuduerel} for historical comments concerning the
next proposition. 

\begin{proposition} \labbel{pu}
\cite{apuduerel} 
Consider the following properties of 
a binary relation $R$.
\begin{enumerate}[1.]
    \item 
$R$ is transitive;
\item
$R$ is reflexive;
\item
$R$ is symmetric;
\item
$R$ is antireflexive, that is, 
$x \mathrel { R } x$ never holds;
\item
$R$ is antisymmetric.
 \end{enumerate}
Then the following statements hold.
  \begin{enumerate}[(A)]
    \item   
For every  $P \subseteq \{ 1,2,3,4,5\} $,
the class $\mathcal K_P$
 of the structures with a binary  relation
$R$ satisfying  the corresponding properties  has  superSAPU.

\item
 For each $P \subseteq \{ 1,2,3,4,5 \} $, let $\mathcal K^f_P$
be the class of structures 
with an added unary operation $f$   
which is $R$-preserving, that is, 
\begin{equation}\labbel{rel}     
\text{ $  x \mathrel { R } y $ implies   
 $  f(x) \mathrel { R } f(y) $.}
 \end{equation}
Then  
$\mathcal K^f_P$
has superSAPU.

\item
 For each $P \subseteq \{ 1,2,3,4,5 \} $, let $\mathcal K^g_P$
be the class of structures 
with an added unary operation $g$   
which is $R$-reversing, that is, 
\begin{equation}\labbel{rev}     
\text{ $  x \mathrel { R } y $ implies   
 $  g(y) \mathrel { R } g(x) $.}
 \end{equation}
Then  
$\mathcal K^g_P$
has superSAPU.

\item
For every  $P \subseteq \{ 1,2,3,4,5 \} $ and
for every pair $F$, $ G$ of sets,
let $\mathcal K^{F,G}_P$
be the class of 
models obtained 
from members of $\mathcal K_P$
 by adding an  
 $F$-indexed set  of unary operations 
 satisfying  
\eqref{rel}
and a
$G$-indexed set of unary operations 
 satisfying  
\eqref{rev}.

Then
$\mathcal K^{F,G}_P$
has superSAPU. 
 \end{enumerate}
  \end{proposition}

Proposition \ref{pu} is proved in \cite{apuduerel}. 
We give here a sketch of the proof, since it
will be used in the sequel.

  \begin{proof}[Sketch of proof]
Assume that 
$\mathbf A$, $\mathbf B$, $\mathbf C$
is a TBA triple (a Triple to Be Amalgamated).  
The proof of (A) is  divided into two  cases.

(A, case a) Suppose that  $ 1 \notin P$.
Then
let $R$ on $A \cup B$ be defined  
by $R= R _{\mathbf A} \cup R _{\mathbf B}$.
Namely, $d \mathrel { R } e $
if and only if 
\begin{equation}\labbel{ei}     
\text{ \underline{either} $d,e \in A$
and   $d \mathrel {  R _{\mathbf A}} e $,
\underline{or} 
 $d,e \in B$
and   $d \mathrel {  R _{\mathbf B}} e $.}
 \end{equation}
Let $\mathbf D = (A \cup B, R)$. 

(A, case b)
Suppose that $ 1 \in P$.
Let $d \mathrel { R } e $
if  either \eqref{ei}, or
\begin{equation} \labbel{al}
\begin{aligned}
&d \in A, e \in B \text{ and there is $c \in C$ such that  }
d \mathrel { R_{\mathbf A} } c \mathrel { R_{\mathbf B} } e, \text{ or }  
\\  
&d \in B, e \in A \text{ and there is $c \in C$ such that  }
d \mathrel { R_{\mathbf B} } c \mathrel { R_{\mathbf A} } e.  
 \end{aligned}  
\end{equation}  
Thus in this case we set
 $R= R _{\mathbf A} \cup R _{\mathbf B} \cup 
(R _{\mathbf A} \circ R _{\mathbf B}) \cup 
(R _{\mathbf B} \circ R _{\mathbf A})$.
Again, let $\mathbf D = (A \cup B, R)$. 

In both cases it can be checked that $\mathbf  D$ 
superamalgamates $\mathbf A$ and $\mathbf  B$ over $\mathbf  C$
and $R$ on $\mathbf  D$ satisfies the required conditions.
This proves  (A).

Given the proof of (A), clauses
(B), (C) and (D) are proved by observing that
 the operations $f$ and $g$  are uniquely defined on $\mathbf  D$, since
$D= A \cup B$, and then by checking that $f$ satisfies
the required properties. 
 \end{proof}

In order to get APU, operations should be unary 
in \ref{pu}(B)-(D). See Proposition  \ref{bin}.

\begin{corollary} \labbel{propj}
The following classes have superSAPU and JEP:

(a) partially ordered sets,
(b) preordered sets,
(c) undirected graphs (sets with a symmetric antireflexive binary relation),
(d) directed graphs (sets with an antireflexive binary relation),
(e)
sets with  an equivalence relation,
(f)
 sets with  a binary transitive relation,
(g)
 sets with  a binary symmetric and reflexive  relation (a \emph{tolerance}).

SuperSAPU and JEP are maintained when
an $F$-indexed family of relation-preserving \eqref{rel}
unary operations and 
a $G$-indexed family of  relation-reversing \eqref{rev}
unary operations are added. 
 \end{corollary}

In Section \ref{pres}  below we shall deal with
AP for structures with many binary relations.
We now present  a theorem 
(from \cite{apuduerel}) which does not seem to follow from the
results in the subsequent sections.

Given two binary relations $R$ and $S$ 
on the same domain,
we say that  $S$
is \emph{coarser} than  $R$ if $R \subseteq S$,
more explicitly, if $a \mathrel { R } b $ implies 
$a \mathrel { S } b $, for all $a$ and $b$ in the domain.
If this is the case, we shall also say that $R$ is \emph{finer}
than $S$. Notice that we shall always  use the expression
``coarser'' in the sense of ``coarser than or equal to''.

\begin{theorem} \labbel{due} \cite{apuduerel}
  \begin{enumerate}[(A)]
    \item   For every pair  $P, Q \subseteq \{ 1,2,3,4,5 \} $,
the class $\mathcal K_{P,Q}$ of structures with 
  \begin{enumerate}
   \item  
a binary relation $R$ satisfying the properties 
from $P$ and
\item 
 a coarser
relation $S$ satisfying the properties 
from $Q$
  \end{enumerate}
 has SAPU.
Actually, superSAPU
holds both with respect to 
$R$ and $S$.

\item 
 SuperSAPU is maintained if we add families of 
  \begin{enumerate}[(i)]
    \item 
unary operations which are both $R$- and $S$-preserving;
\item
unary operations which are both $R$- and $S$-reversing;
\item
unary operations which are  $R$-preserving;
\item
unary operations which are $R$-reversing;
  \end{enumerate} 

SuperSAPU is maintained if we
consider the subclass consisting of those structures
satisfying any set of properties among those 
we shall list in Propositions \ref{un}, \ref{agg}  
and Theorem \ref{gen}   below.
\item
 On the other hand, 
the class of structures with a transitive relation $R$,
a coarser binary relation $S$ and an $S$-preserving function $f$  
has not AP (here we do not include the condition that $f$ 
is  $R$-preserving).

\item
Similarly, the class of structures with a partial order $\leq$,
a coarser symmetric and reflexive relation $S_1$ and an 
$S_1$-preserving function $f$  
has not AP.
\end{enumerate} 
 \end{theorem}

\begin{prop} \labbel{div+} \cite{apuduerel}
The following theories do not have  AP.
  \begin{enumerate}[(a)]
    \item   
The theory of an antisymmetric relation $S$ 
with two partial orders $\leq$ and $\leq'$
both finer than $S$. 

\item
 The theory of an antisymmetric relation $S$ 
with two transitive relations
both finer than $S$.

\item
 The theory of a partial order with
a coarser antisymmetric relation $S$
and a bijective $S$-preserving 
unary operation.  
\end{enumerate} 
 \end{prop}   
 
Notice that either (a) or (b) in Proposition \ref{div+} 
shows that a universal Horn theory in a pure relational language
does not necessarily have AP:

We now recall some results 
from \cite{lop}.
A generalization of Corollary \ref{propj}(2)
to linear orders holds, but, rather unexpectedly,
it holds only in the case of  just one additional operation.

If $f$ is an unary operation on a poset,
we say that an element $c$ is a \emph{center},
or a \emph{fixed point} 
for $f$ if $f(c)=c$.
In general, when we refer to a \emph{center}
$c$ without further specifications, we shall
mean that   $c$ is a center for all the operations
under consideration.

\begin{theorem} \labbel{lopp}
\cite{lop}
The following classes have SAPU.
  \begin{enumerate} 
   \item 
 The class 
 of linearly ordered sets with one order preserving
 unary operation.

\item
The class of linearly ordered sets with 
one order reversing unary operation with a center.

\item
For every set $ F$, the class 
 of 
linearly ordered sets with an $F$-indexed
family of  order automorphisms.

\item
For every pair $F$ and $G$ of sets, the class of
linearly ordered sets 
with an $F$-indexed family of  order
automorphisms and a $G$-indexed family
of bijective order reversing unary operations,
all operations from both families with a common center.
\end{enumerate} 

The following classes have APU
but not SAP.

  \begin{enumerate}
\setcounter{enumi}{4}

\item The class 
of linearly ordered sets with one strict order preserving
unary operation.

\item
 The classes  
 of linearly ordered sets 
with one order reversing, resp.,
one strict order reversing unary operation.

  \end{enumerate}

The following classes have not AP.

  \begin{enumerate}
\setcounter{enumi}{6}

\item
 The classes  
 of linearly ordered sets 
with two order preserving, resp.,
two strict order preserving unary operations.

\item
More generally, the classes 
of linearly ordered sets with 
 two 
unary operations and
with each operation  
 either
order preserving, or
strict order preserving, or
order reversing, or 
strict order reversing.
  \end{enumerate}
 \end{theorem}

Of course, when dealing with
linearly ordered sets,  we cannot have superSAP.
Given $\mathbf A$, $\mathbf B$ and  $\mathbf C$
linearly ordered sets to be amalgamated, if $a \in A \setminus C$,
 $b \in B \setminus C$ and 
$\{ c \in C \,  \mid  c < _{ \mathbf A} a\,\} = 
\{ c \in C \,  \mid  c < _{ \mathbf B} b\,\}$,
then in any strong amalgamating \emph{linear}
order we have either $a<b$ or $b<a$,
but no such relation is witnessed by means of some $c \in C$.

Theorem \ref{lopp}(1) does not hold for binary 
operations which are order preserving  
on each component, as we are going to show in the next example.

\begin{example} \labbel{exbin}
The class of linearly ordered sets with 
a binary 
operation which is order preserving  
on each component has not AP.

Take $C= \{ c \} $,
$A= \{  a,c\} $,
with $c < a$ and  a binary function $f$  defined by
$f(a,c)=c$,  $f(c,a)=f(a,a)=a$, and 
$B= \{  b,c\} $,
with $c<b$ and  $f$  defined by
$f(c,b)=c$,  $f(b,c)=f(b,b)=b$.

In any amalgamating algebra 
we cannot have $b \leq a$, since then 
$b=f(b,c) \leq f(a,c)=c$, as 
 $f$ is requested to  be order preserving
on the first component,
Symmetrically, we cannot have $a \leq b$, 
hence $\mathbf A$, $\mathbf  B$ and $\mathbf  C$ 
cannot be amalgamated to a linear order.
\end{example}

\section{Preservation conditions.} \labbel{pres}

We now show 
how to construct new theories with SAPU
 starting from some theories with the property.
For theories in a purely relational language,
Part (a) in the following proposition is folklore.
Then Part (b) is an immediate consequence, but it seems 
to have not received due attention in the literature.

\begin{proposition} \labbel{T}
(a)
If $( T_i) _{i \in I} $ 
is a sequence of theories in pairwise disjoint languages
and each $T_i$ has SAPU, then 
$ T=  \bigcup_{i \in I} T_i $ has SAPU.

(b)
If $T_1$ is a theory in some language $\mathscr L$ 
and $T_1$ has  (super)SAPU, then $T_1$ has (super)SAPU even when 
considered as a theory in some language $\mathscr L' \supseteq \mathscr L$.    
\end{proposition}  

\begin{proof}
(a) As in Remark \ref{iso}, suppose that   $\mathbf A$, $\mathbf B$, $\mathbf C$
is a TBA triple  consisting of models of $T$.  For each $i \in I$,
the reducts to the language of $T_i$ can be amalgamated 
to a model over $A \cup B$.
Since the languages are pairwise disjoint, 
we get a model of the whole $T$ over $A \cup B$.

To prove (b), let $T_2$ 
be the empty theory in the language $\mathscr L' \setminus  \mathscr L$.
By Observation \ref{obs}(c), $T_2$ has SAPU.
Then apply the first statement.
If $T_1$ has superSAPU, then superSAPU is maintained by construction. 
 \end{proof} 

See Proposition \ref{usup} below for a result analogue to 
Proposition \ref{T}(a) for the superamalgamation property.

We shall see in Proposition \ref{senzau}
that the analogue of Proposition \ref{T} holds when we replace  SAPU  with
SAP. However, the proof of Proposition \ref{T}
is much simpler and the method of proof can be applied to 
more situations, see Theorems \ref{mul} and \ref{gen}.  
Moreover, there are results holding for SAPU
but not for (S)AP: compare Proposition \ref{un} with  
 Example \ref{uneed} below.
See also Remark \ref{card} and Example  \ref{***}.
The assumption that the $T_i$'s have SAPU 
in Proposition \ref{T} cannot be weakened to APU.  
See Example \ref{necess} below. 
Compare also Example \ref{uneed}(b). 
A slightly more general version of 
Proposition \ref{T} is stated as Proposition \ref{L} 
below.

In order to present the following results
in due generality, we need to introduce some
terminology and notation.
An \emph{($I$-indexed) multiposet}
is a set endowed with a family $( \leq_i) _{i \in I} $ of partial orders.
It is immediate from Proposition \ref{T} and Corollary \ref{propj}
that, for every set $I$, the class of all the    
$I$-indexed multiposets has SAPU.
We are going to prove a more general fact 
about multiposets 
on which some coarseness conditions
are assumed.
In contrast with Theorem  \ref{due},
here we impose no bound on the cardinality of $I$. 
A \emph{coarseness condition} on an
($I$-indexed) multiposet is 
a condition of the form ``$ \leq_i$ is coarser
than  $ \leq_j$'', for some pair
$(i,j)$, with  $i, j \in I$.  
Thus a family of coarseness conditions 
is (represented by) a subset 
$\mathcal F$  of $I \times I$: we are asking that 
$ \leq_i$ is coarser
than  $ \leq_j$ for all pairs 
 $(i,j) \in \mathcal F$.

If $\mathcal K$ is a class of structures,
let us denote by $\mathcal K^{fin}$ 
the class of the finite members of $\mathcal K$.

\begin{theorem} \labbel{mul}    
  \begin{enumerate}[(1)] 
   \item 
 For any index set $I$, the class of all
$I$-indexed multiposets satisfying any given 
 family of  coarseness conditions on $I$ has SAPU,
actually, superSAPU with respect to each $\leq_i$.

\item
Suppose that  $J \subseteq I$,
$\mathcal F$ is a family of coarseness
conditions on $I$ and
let $\mathcal K$ be the class of all
$I$-indexed multiposets which satisfy the
  coarseness conditions in $\mathcal F$ 
and  such that all orders $\leq_j$ with $j \in J$  
are linear. Then $\mathcal K$ 
 has SAPU.
\item
In particular, the theory of 
 a partial order $\leq$  together with a
linearization of $\leq$ has  SAPU.
 \end{enumerate} 

If $I$ is finite, then in each case
the class of finite structures has a Fra\"\i ss\'e limit 
 $\mathbf M$ 
and  the
first-order theory $Th(\mathbf M)$ is
$ \omega$-categorical and has quantifier elimination; moreover, 
$Th(\mathbf M)$
is the 
model completion of the
first-order theory axiomatizing the class
under consideration.
\end{theorem}

\begin{proof} 
(1) Formally, the theorem   is not a consequence 
of Proposition \ref{T}.  
However, if we apply the proof of 
Proposition \ref{pu} simultaneously for all the relations
involved,
we get a structure  for the appropriate language.
Since the order relations are all
 transitive, we are always in case b, 
hence coarseness is preserved.

(2) If $\leq_j$ is a linear order
on  $\mathbf A $, $  \mathbf B $ and $  \mathbf C $, then 
the proof of Proposition \ref{pu}
generally provides only a partial order $\leq_{j,\mathbf D}$;
however, any linearization of $\leq_{j,\mathbf D}$
works so as to get an amalgamating structure
with a linear order.
If some coarseness condition asserts that 
$\leq_j$ is coarser than $\leq_i$,
and $i \in I \setminus J$,
that is, $\leq_i$ is assumed to be partial,   then,
as in (1),
$\leq_{j,\mathbf D}$ is coarser than
$\leq_{i,\mathbf D}$ and thus any linearization
of $\leq_{j,\mathbf D}$ is coarser than
$\leq_{i,\mathbf D}$.
The other case is trivial:
if some coarseness condition asserts that 
$\leq_i$ is coarser than some linear order  $\leq_j$,
then necessarily ${\leq_i} = {\leq_j}$,
hence we can take  ${\leq_i} = {\leq_j}$
in $\mathbf  D$.

To prove the last statement,
first notice that, for every class $\mathcal K$ under consideration,
the class $\mathcal K^{fin}$ has AP, 
since we can amalgamate into union. JEP follows since
here we are allowed to consider an empty $\mathbf  C$.
Then use \cite[Theorems 7.1.2 and 7.4.1]{H}. 
 The finiteness assumption  is necessary
in order to have a countable number of structures under
isomorphism.
\end{proof}

We now turn to another method
which produces theories with SAPU,
starting from theories satisfying the property.
The proof is trivial, but the method is useful.
For short, SAPU is preserved if we only consider
those models which satisfy some given set of universal-existential sentences
in which only one variable is bounded by the universal quantifier.

\begin{proposition} \labbel{un}
Fix some language for all the sentences and the models under consideration.

Suppose that  $\Sigma = \{  \sigma_i \mid i \in I \}  $ is a set of 
universal-existential sentences
in which at most one variable is bounded by the universal quantifier, namely,
 sentences
of the form
\begin{equation} \labbel{sent}  
\forall x \exists y_1 y_2 \dots  \ \varphi _i
\qquad \text{\rm or } \qquad
 \forall x \varphi_i
\qquad \text{\rm or } \qquad
\exists y_1 y_2 \dots \varphi _i
  \end{equation}     
where in  each case  $\varphi_i$  is
quantifier-free.

(a) If $\mathcal K$ is a class of structures with
(super)(S)APU, then the class $\mathcal K'$
of all structures in $\mathcal K$ which satisfy 
$\Sigma$ has (super)(S)APU.

(b) In particular, if $T$ is a theory with (super)(S)APU
 then $T \cup \Sigma $ has (super)(S)APU. 

(c) Any theory with only axioms of the form
\eqref{sent} has  SAPU, superSAPU, in the presence of a binary relation.
\end{proposition}  

\begin{proof} 
 If $ D =  A \cup  B$  
and both  $\mathbf A $ and $  \mathbf B$ satisfy
some sentence of the form \eqref{sent}, then
$\mathbf D$  satisfies such a sentence, if 
both $\mathbf A$ and $\mathbf  B$ are substructures of $\mathbf  D$.
Notice that at most one variable 
is bounded by $\forall$. 

 The last statement follows from Observation \ref{obs}(c).
\end{proof}

For example, by a formula of the form \eqref{sent}
we can express the condition that some unary operation $f$
satisfies identically $f(f(x))=x$, or $f(x) \geq x$.
We can say that two unary  operations $f$ and $g$
are comparable, $f(x) \geq g(x)$. We can say that some function
is surjective, or even that some relation is surjective
with respect to 
some component, e.~g., 
$\forall x_1\exists x_2, \dots x_n \, R(x_1,x_2, \dots x_n )$, etc.

\begin{example} \labbel{kisi} 
A.  Kisielewicz \cite{Ki} presents examples 
of varieties without nontrivial finite algebras.
A simple example is the variety $\mathcal V$  with three unary operations
$f$, $g$ and  $h$ satisfying
$fgh(x)=x$ and $fh(x)=fh(y)$ identically.
 Indeed,   the 
first   identity implies that  $f$  is  surjective  and  $h$
  is  injective. The  second  identity 
implies that either   $h$  is  not  surjective or $f$  is constant.  

As an application of Proposition \ref{un} 
we show that $\mathcal V$ has SAPU (for nonempty
algebras). Formally, Proposition \ref{un}
does not apply to the the second identity; however, 
if we introduce a new constant $c$ and we replace the second identity
by $fh(x)=c$, then, for nonempty algebras, we get exactly
the same morphisms and embeddings, hence Proposition \ref{un}    
can be applied.
\end{example}

The ``U'' in SAPU is necessary in Proposition \ref{un},
see Example \ref{uneed} below. 
In Proposition \ref{un}  it is necessary to assume that
in \eqref{sent} at most one variable    
is bounded by $\forall$. See Remark \ref{inf}(d) or  Example \ref{1} below.

\begin{remarks} \labbel{inf}    
(a) We do not need the sentences in 
Proposition \ref{un} to be finitary,
they might possibly  be infinitary.
We might have infinitely many variables $y_j$,
as far as \emph{at most one} variable    
is bounded by $\forall$ and the $\varphi_i$'s  
are quantifier-free.

(b) There are first-order sentences
for which the statement of Proposition \ref{un}
holds (limited to SAPU), but which do not have the form \eqref{sent}.
For example, if some class $\mathcal K$ in the language
with a unary operation $f$ has SAPU, then the subclass
of those structures in $\mathcal K$ in which
$f$ is bijective has still SAPU.
We shall show in Example \ref{injbij}
that APU is not preserved by adding a sentence saying
that some function is bijective; actually, 
AP might be destroyed.

(c) The assertion that $f$ is bijective     
cannot be expressed by a sentence of the form \eqref{sent},
however, it can be expressed by a particularly
simple second-order sentence, since 
$f$ is bijective if and only if $f$ has   an inverse. 
Hence the solution to the problem of finding the most general
form of Proposition \ref{un}
(see Problem \ref{probun} below) 
might involve second-order sentences.

(d) Notice that  the assertion that $f$ is surjective     
can be actually expressed by a sentence of the form \eqref{sent}.
On the other hand, SAPU for some class is not preserved
by adding the condition that some operation is injective.   
Adding such a condition might even
destroy AP: 

For example, consider the class $\mathcal K$ of all structures
with a unary operation $f$ and a unary predicate $V$.
The class $\mathcal K$ has SAPU by 
Observation \ref{obs}(c). 
Let $\mathbf  C$ be $\mathbb N$ 
with $f$ the successor operation.
Let $A= \mathbb N \cup \{ a \} $ with
$f(a)=0$    
and $B= \mathbb N \cup \{ b\} $
with $f(b)=0$.
If $f$ is still to be injective, we cannot have SAP.
If we further set $V(a)$ and  $ \neg V(b)$,
even AP fails.    

Hence SAPU is not preserved by adding the condition that 
some operation is injective.
\end{remarks}

\begin{problem} \labbel{probun}
(a) Characterize those sets of sentences $\Sigma$ 
for which the analogue of Proposition \ref{un} holds,
either for SAPU or for APU. Notice that the two cases are distinct, by
Remark \ref{inf}(b) and   Example \ref{injbij}.

(b) Are there even more sentences for which
Proposition \ref{un}(c) holds? Some affirmative answers are 
provided in \cite{apuvar}.

(c) More generally, characterize those sets of
(not necessarily first-order)
 properties $H$
such that, whenever $\mathcal K$ is a class with (S)APU,
then the subclass of $\mathcal K$ consisting of those structures 
satisfying the properties in $H$ has still (S)APU.
See Proposition \ref{agg}  
 for more examples of such properties.

(d) Solve the above problems for (S)AP in place of (S)APU.
From Example \ref{uneed} below,
we see that Proposition \ref{un}, as stated, fails
for SAP  in place of SAPU.
However, it might happen that some version of \ref{un}
holds for SAP, when a more restricted set of sentences
is taken into account.
There are trivial cases,
for example,  Proposition \ref{un}
holds for AP, when we restrict 
\eqref{sent} 
to existential sentences.

As  a test case, is it true that if $T$ has SAP,
then the theory which further asserts that some unary function is
bijective has still SAP?
 \end{problem}

Recall that 
a \emph{closure operation}
on some poset $\mathbf A$ 
is an order preserving unary operation 
 $f$ such that 
$x \leq f(x)=f(f(x))$, for every $x$.
See \cite{E} for further information
and pictures. 
In the presence
of a semilattice operation, some authors include 
an additivity  requirement
in the definitions of  closure.
We shall adopt the more general convention
\cite{E} according to which no additivity
assumption is made.

An (\emph{antitone}) \emph{involution}
is an (order-reversing) unary operation $'$
such that   $x''=x$.

\begin{corollary} \labbel{clp}
  \begin{enumerate}[(A)]    
\item  
For every pair $F$, $G$ of sets, the class of posets 
 with an $F$-indexed family of closure operations
and a $G$-indexed family of antitone involutions
 has
superSAPU. 

\item
The class of 
linearly ordered sets with one closure operation 
has SAPU.
\item
The class of 
linearly ordered sets with  one antitone involution 
has APU.

\item
The class of  linearly ordered sets with two  closure operations has
not AP.

\item
The class of 
 linearly ordered sets with two antitone  involutions  has
not AP. 

\item
The class of
 linearly ordered sets with a family of antitone involutions
with a common fixed point
 has
SAPU.
 \end{enumerate}
 \end{corollary} 

\begin{proof}
(A) follows from Proposition \ref{pu}(D)
 and 
   Proposition \ref{un},
 considering  sentences of the form
$\forall x \ x \leq f(x)=f(f(x))$
and $\forall x \ x''=x$.

(B) and (C)  follow from 
   Proposition \ref{un} and 
\cite[Theorems 3.1(a) and 4.3(a)]{lop},
which have been  recalled in
Theorem \ref{lopp}. 

(D) appears in \cite[Remarks 3.2]{lop}.

(E) Case (b)(iii)
in the proof of \cite[Theorem 4.3]{lop} 
provides a counterexample, though
involutions are not explicitly mentioned
in \cite{lop}.

(F)  follows from
\cite[Theorem 4.5]{lop}, 
using again  Proposition \ref{un}.
Notice that the assumption that two involutions $'$
and $^*$ have a common   fixed point can be expressed by
the sentence  $\exists x \ (x'=x\  \&\  x^*=x)$,
which  has the form \eqref{sent}.
 \end{proof}

There are also many non first-order properties which preserve SAPU.
We present some examples. 

If $R$ is a binary relation on some set $A$, an \emph{$R$-antichain}
is a subset $X$ of $A$ such that not  $a \mathrel { R } b $,
for every pair of distinct elements $a , b \in X$.   

It is trivial that if $\mathbf A$, $\mathbf  B$, $\mathbf  C$ 
is a TBA triple of connected graphs, $\mathbf  C$ is nonempty
and $\mathbf  D$ is an amalgamating structure, then $\mathbf  D$ is
connected, too. The notion of connectedness
can be generalized in various ways in model theory.
We present a quite general version.

\begin{definition} \labbel{conne}    
If $\mathbf A$ is a model, two elements 
$a, b \in A$ are \emph{adjacent} 
if $R(a_1, a_2, \dots, a, \dots, b, \dots)$
holds, for some relation $R$ 
in the language of $\mathbf A$ and some
$a_1, a_2, \dots \in A$.
We are not assuming that $a$ occurs before  $b$
in the expression $R(a_1, a_2, \dots, a, \dots, b, \dots)$.  
If the above holds, we also say that 
$a, b \in A$ are \emph{$R$-adjacent} 
and also that there is an 
\emph{$R$-$m$-$n$-directed edge from $a$ to $b$},
where $a$ occurs in the $m$th position
and  $b$ occurs in the $n$th position
in $R(a_1, a_2, \dots, a, \dots, \allowbreak  b, \dots)$.

A structure $\mathbf A$ is \emph{connected}
if every pair of elements of $A$ can be connected by a path
consisting of adjacent elements.
In other words,  $\mathbf A$ is \emph{connected}
if the transitive closure of the adjacency
relation on $A$ is the
largest relation $A \times A$. 
The structure $\mathbf A$ is \emph{$R$-connected}
if every pair of elements of $A$ can be connected by a path
consisting of $R$-adjacent elements.
 
Still more generally, let $\mathcal R$
be a family of triples of the form
$(R,m,n)$, where $R$ varies among the relations in the language
and $m,n$ $\leq$ the number of arguments of $R$. 
An $\mathcal R$-path
is a sequence $a_1, \dots, a_h$
such that, for every $i < h$,
there are some
$(R,m,n) \in \mathcal R$ 
and an $R$-$m$-$n$-directed edge from $a_i$ to $a _{i+1}$.
A model $\mathbf A$  is \emph{$\mathcal R$-connected} 
if, for every $a,b \in A$,
there is an $\mathcal R$-path
with initial point  $a$ and final point $b$. 
See \cite{Co} for related notions.

For the purpose of the above definitions 
of connectedness, we can take into account also
function symbols: think of an $n$-ary function
$f$ as an $n{+}1$-ary relation given by 
$R(a_1, a_2, \dots, a_n, a _{n+1} )$
if    $f(a_1, a_2, \dots, a_n)= a _{n+1} $.
\end{definition}

In the statement of the 
following proposition  $\lambda$ is an infinite cardinal.

\begin{proposition} \labbel{agg}
If $\mathcal K$ is a class of structures with (S)APU,
then, for any set of properties chosen from the  list below,
the subclass of $\mathcal K$ consisting of those structures in $\mathcal K$
which satisfy the chosen properties has (S)APU.
  \begin{enumerate}[({A}1)]    
\item
The domain is finite (or has cardinality $< \lambda$).
\item
The domain is finitely generated  (generated by a set of  cardinality $< \lambda$).
\item
For some binary relation  $R$ 
assumed to be a partial order:
$R$ is well-founded
(has no strict descending chain
of length $ \geq \lambda$) (has no strict ascending chain
of length $\geq \lambda$).
\item 
For some binary relation $R$, there is no
infinite 
$R$-antichain
(there is no 
$R$-antichain
of cardinality $< \lambda $)
 
\item 
(only for SAPU)  
For some natural numbers $n, m$
fixed in advance, the domain has cardinality 
$kn+m$, for some $k$.  
\end{enumerate}

As far as the following properties are concerned, we 
consider the version of (S)APU
in which the bottom structure $\mathbf  C$ is assumed to be nonempty
(compare the second paragraph in Section \ref{sec}).
 \begin{enumerate}[({C})]
  \item 
The structure is connected
($R$-connected, for some relation symbol $R$)
($\mathcal R$-connected, for some family of triples $\mathcal R$ 
as in Definition \ref{conne}).
  \end{enumerate} 

More generally, if $V$ is  a unary predicate,
we can consider anyone of the above properties when
restricted to the domain $\{ \, x  \mid V(x) \,\}$
of $V$. When applying condition (C) 
we should assume that the domain of $V$
is nonempty. 
 \end{proposition}

The proof of Proposition \ref{agg} is immediate;
in fact, if some of these properties hold
in $\mathbf A$, $\mathbf B$ and  $\mathbf C$,   
then the property holds in some amalgamating structure, 
since we can construct it
on $A \cup B$. 
However, the corresponding statements
generally fail when SAPU is weakened to SAP.
For instance, see Example \ref{1a}(c).  

Recall that a 
 \emph{well partial order},  or \emph{wpo}
is a  well-founded partial order
 without infinite antichains.
Hence, because of Proposition \ref{agg}(A3)(A4), 
if some class with (S)APU has a partial order
relation $\leq$, then the subclass of those structures in which
 $\leq$ is a wpo has still (S)APU.

Some assumption on the relation $R$ 
in Proposition \ref{agg}(A3) is necessary: 
see Example \ref{exch}.

Theorem \ref{mul} and Corollaries \ref{propj} and  \ref{clp}
have not the most general form, rather, they are 
just exemplifications. We are now going to state  a very general result
which can be obtained by 
the present methods.

By the proof of  Proposition \ref{pu},
we can deal with structures with a family of transitive relations, each
relation satisfying any set of properties chosen from 2.- 5.
Structures with partial orders are just an instance of this more general case. 
In particular, we can deal with  any class of structures with many preordered sets,
 many equivalence relations, or even  simultaneously
 posets, 
preordered sets
and equivalence relations.

Even in this general setting 
we can define a \emph{coarseness condition},
namely,
a condition of the form $R_j \subseteq R_i$.
That is, if $x \mathrel { R_j } y$,
then $x \mathrel { R_i } y$.

Due to Propositions \ref{pu}(D) and \ref{T}(b),
we can add a family of unary operations,
as well as conditions asking 
that some operation is 
$R$-preserving or $R$-reversing.
We can add sentences of the form
\eqref{sent}, in view of Proposition \ref{un}.    
In conclusion,  here is a general result we have got.

\begin{theorem} \labbel{gen}
All the classes $\mathcal K$  described below 
 have SAPU; actually, superSAPU with respect to 
all the binary relations $R_i$ involved.

We assume  that $\mathcal K$ is the class of models
for some theory $T$ 
in a language with
a sequence $( R_i) _{i \in I} $
of binary relation symbols and  
a sequence $( f_h) _{h \in H} $
of unary function symbols.
We require that  $T$ asserts that each 
$R_i$ is a transitive relation. Moreover, 
$T$ is allowed to contain some axioms, possibly none,
from the list below
(each axiom might appear for as many indices as wanted).
  \begin{enumerate}   
 \item 
Some relation $R_i$ is reflexive;
\item 
Some relation $R_i$ is symmetric;
\item 
Some relation $R_i$ is antireflexive;
\item 
Some relation $R_i$ is antisymmetric;
\item
Some $f_h$ preserves some $R_i$;
\item
Some $f_h$ is $R_i$-reversing;
\item
Some $R_i$ is coarser than some $R_j$;
\item
Any sentence of the form \eqref{sent}, in particular,
the universal closures of sentences of the form
  \begin{enumerate}    
\item 
$ f_h(x) \mathrel { R_i }   f_k(x)$, 
\item
$ x \mathrel { R_i }   f_h(x)$, 
\item
$ f_h(f_h(x)) = f_h(x)$, 
\item
$ f_h(f_k(x)) = f_k(f_h(x))$,
etc.
 \end{enumerate}   
\item
Some $f_h$ strictly preserves some $R_i$,
provided $T$ asserts that $R_i$ is a partial order;
\item
Some $f_h$ is strictly $R_i$-reversing,
provided $T$ asserts that $R_i$ is a partial order.
\end{enumerate}

Furthermore, we are allowed to expand the language by adding any number
of  symbols of any kind, as far as the
axioms involving such new symbols are only of the form
\eqref{sent}.   

Still more generally, for each class $\mathcal K$ as above, the subclass
of those substructures satisfying any given set of conditions
taken from Proposition \ref{agg} has SAPU. 
 \end{theorem} 

\begin{proof}
As in the proof of \ref{mul},
apply the proof of 
Proposition \ref{pu} individually for each
relation and function, joining everything in the model 
$\mathbf D$. By the proof of Proposition \ref{pu}, if any
sentence equivalent to some condition from (1) - (6) is included in $T$,
then $\mathbf D$ satisfies this sentence.
 By assumption, all the relations
are transitive, hence we are always in case b
in the proof of \ref{pu}; this implies that
coarseness is preserved. This argument takes care of (7).
Clause (8) follows from Proposition \ref{un}.
 Clauses (9) and (10) follow  from Remark 2.2 in \cite{apuduerel},
where we check that the construction in case b
in the proof of Proposition \ref{pu}  commutes in passing from some partial
order to the corresponding strict order.

As far as the penultimate statement is concerned,
use Proposition \ref{T}(b)
and  again Proposition \ref{un}. The last statement is immediate from
Proposition \ref{agg}.
 \end{proof}

\begin{remark} \labbel{fragen}
In most cases, the class $\mathcal K^{fin}$
has a Fra\"\i ss\'e limit, for a class $\mathcal K$ as considered
in Theorem \ref{gen}.
Compare  Theorems  \ref{mul}.

However, there are some limitations.
First, we generally need the language to be finite, 
in order to get only a countable number of finite models
modulo isomorphism.
Second, only universal sentences of the form \eqref{sent}
can be considered, if we want  hereditariness \cite[7 (1.1)]{H}   to be preserved.  
 Finally, some of the conditions mentioned
in Theorem \ref{gen} might destroy JEP,
for example, adding constants to the language, or
using Clause (C) from Proposition \ref{agg}.  
However, 
we retain JEP when restricted to any equivalence class
as described in Remark \ref{jepeq},
hence we have a Fra\"\i ss\'e limit for each equivalence class.

We leave the details to the interested reader.
 \end{remark}

\begin{remark} \labbel{bup}
If we merge transitive relations 
with relations which are not
supposed to be transitive, then 
coarseness is not always preserved
when trying to amalgamate structures;
see Proposition \ref{div+}. 

In particular,  Theorem \ref{gen}
does not necessarily hold when transitive relations 
are merged with nontransitive relations,
if  we ask for coarseness conditions as in (7).
The case of just two relations, as given by
Theorem \ref{due},  is a notable exception.
Of course, as in the first lines of the 
proof of \ref{due}, if we compare
any number of relations which are not required to be transitive,
then coarseness is preserved.
 Similarly,   we can ask that some relation is finer than
another \emph{transitive} relation, since the definition of $R$ in the proof of 
Proposition \ref{pu}, case b,  
always produces a coarser relation, in comparison
with $R$ as defined in case a. 

Quite unexpectedly, a common generalization
of  Theorems \ref{due} and \ref{gen}(1)-(4), (7)-(8)
holds 
even when dealing with relations supposed to be 
coarser than other transitive relations.
We shall show that the only obstacle to 
amalgamation are antisymmetric relations
supposed to be  coarser than a pair of incomparable transitive relations;
essentially, the only counterexamples to amalgamation
are given by the examples we have described
 in Proposition \ref{div+}.  
However, in the general case which unifies 
Theorems \ref{due} and \ref{gen},
  the possibility of adding relation-preserving operations
is not as neat as in Proposition \ref{pu}(D) or Theorem 
\ref{gen}(6), (7), (9), (10).
Samples witnessing this difficulty
are given here in Theorem \ref{due}(C)(D)
and Proposition \ref{div+}(c). 
We shall present further details elsewhere.
\end{remark}

\begin{remark} \labbel{pow}    
Theorem \ref{gen} is quite powerful.
We present a simple example.

A \emph{bounded poset} is a poset with a maximum and
a minimum element, both elements interpreted as constants.
We could repeat all the above arguments getting corresponding
theorems for bounded posets.

However, the results follow automatically
from Proposition \ref{un} and Theorem \ref{gen}.
Indeed, by Propositions \ref{pu} and  \ref{T}(b), the class of posets
in the language with two constant symbols added has SAPU. 
Then the assertion that, say the constant $1$ is interpreted as a maximum
can be expressed by the sentence 
$\forall x \ x \leq 1$, having the form  \eqref{sent}
from Proposition \ref{un}.

The above arguments apply in the same way 
in order to show that bounded posets
with families of (strict) order preserving (or reversing)
unary operations have SAPU.
Of course, we can also assume only the existence of a maximum, or only
the existence of a minimum (provided, as above, each one is
interpreted as a constant).
\end{remark}

\section{Binary and $n$-ary operations.} \labbel{binbin}

Reading the previous sections, the reader
might expect that (S)APU
 is a phenomenon typical of relational structures only,
 with  possibly unary operations added.
By Observation \ref{obs}(c) and Proposition \ref{T}(b),
we may have SAPU 
when no axiom involves binary or $n$-ary operations, 
but, in a sense, this is a trivial case. 
However, there are many examples of classes
with  $n$-ary operations, for  $n \geq 2$,
and sharing SAPU.

Some varieties of groupoids with SAPU
are obtained in \cite[Proposition 5.19]{Ke}. 

\subsubsection*{Varieties defined by linear equations.} \labbel{v} 
We now present a general result
from \cite{apuvar} dealing with varieties. 
Recall that a variety $\mathcal V$  is a class of 
nonempty structures
for a language without relation symbols,
such that $\mathcal V$  can be defined by equations, i.e., 
universal closures of atomic formulae.
An equation is \emph{linear} 
if there is at most one occurrence of  an operation
symbol on each side (constants are not counted
as operations, here). Notice that the terminology
is not uniform in the literature. Examples of linear equations
are  $f(x,x,y)=y$, $f(x,x,y)= g(x,x,y)$,
 $h(x,c)=x$ or $f(c,y,y,z)= g(x,x,d,y)$. 

 Linear equations
are important because they
 are almost invariably encountered in the definition of Maltsev conditions,
for examples, the conditions characterizing congruence permutability,
distributivity, modularity\dots \ 
See, e.~g.,  \cite{Be}  for details.
See \cite{apuvar} for further comments, examples
and related results.

\begin{theorem} \labbel{apuvarthm}
\cite{apuvar} Any variety which can be defined by a set of linear equations
has SAPU.
 \end{theorem} 

\begin{proof} (Sketch)
Fix some arbitrary element
$ \bar{d}$  of $D=A \cup B$.
For every operation $f $
and $d_{1},d_{2}, \dots \in D$,
set 
$f(d_{1},d_{2}, \dots  ) = \bar{d}$,
unless the value of 
$f(d_{1},d_{2}, \dots  )$
is forced by some identity to be satisfied,
or by the requirement that 
both $\mathbf A$ 
and $\mathbf  B$ should embed into $\mathbf  D$.
Check that such conditions do not clash
and that if $f(d_{1},d_{2}, \dots  ) =
g(e_{1},e_{2}, \dots  )$
is (the evaluation of) an identity to be satisfied,
then the value of 
$f(d_{1},d_{2}, \dots  ) $ is not forced
if and only if the value of  $g(e_{1},e_{2}, \dots  )$
is not forced (this might be cumbersome, in general).
If this is the case,
 both values are equal to $ \bar{d}$,
hence the identity is satisfied. 
Full details appear in \cite{apuvar}. 
\end{proof}

We now  present  other examples
of  structures with a binary operation and with SAPU.

\subsubsection*{Directoids and related structures.} \labbel{dis}
A \emph{directoid}
is a set with a binary operation 
$\sqcup $ 
such that 
the identities 
$x \sqcup x = x$, 
$x \sqcup y = y \sqcup x $
and 
$x \sqcup ((x \sqcup y) \sqcup z) = 
(x \sqcup y) \sqcup z$ hold.
Directoids are an algebraization of directed sets,
see \cite{CL,CGGKLP} for details.
The above directoids are the ``commutative'' ones;
there is also a noncommutative version,
but notice that the terminology in the literature is not uniform.
Most of the following results hold also in the
noncommutative case.  

A \emph{bidirectoid} 
is a set
with two  directoid operations satisfying  the  absorption laws
$a \sqcap (a \sqcup b)= a$
and  $a \sqcup (a \sqcap b)= a$.
Bidirectoids corresponds to 
posets which are both upward and downward directed.

A \emph{maximum}
for a directoid is an element $1$
such that   $x \sqcup 1=1$.
A \emph{minimum}
is an element $0$
such that $x \sqcup 0=x$.    
When we speak of a  \emph{directoid
with a maximum (minimum)},
we simply assert that such a maximum
(minimum) exists, but we do not assume that it is interpreted
by a constant.
In particular, embeddings of directoids with a maximum
need not preserve maxima. 
   
On the other hand,
an \emph{upper bounded directoid (bounded directoid)}
is a directoid with a maximum (and a minimum)
interpreted as constant(s).
In this case embeddings are supposed to preserve maxima
(and minima).

The  following theorem 
is immediate from the proof of 
\cite[Theorems 10 and 11]{CGGKLP}.
The case of bounded bidirectoids is proved in essentially the same way.
See \cite{CGGKLP} for the definition of an \emph{involutive}
directoid. 

\begin{theorem}  \labbel{cgk} 
\cite{CGGKLP}
The classes of upper bounded directoids, 
bounded directoids, bounded bidirectoids and of 
bounded involutive directoids 
have SAPU.
\end{theorem}  

In \cite{CGGKLP}
it is also proved that the classes of (not bounded) directoids and bidirectoids
have SAP. The proof does not give SAPU, however, since 
a new element is added to the union of the amalgamating structures.
We shall prove in \cite{apudir} that this new element
is necessary in general, but we can do without adding new elements
in the finite case. 

As above, when we speak of a 
posets with a maximum,  the maximum is not interpreted
as a constant; in particular, embeddings need not preserve maxima.
The case when maxima are interpreted as constants has been dealt
with in Remark \ref{pow}. 
Similarly, a finite directoid has necessarily a maximum, 
but we do not require that embeddings preserve maxima.
As we mentioned, maxima are required to be preserved only in the class of
upper bounded directoids.
See \cite{apudir} for the proof of the next proposition.

\begin{proposition} \labbel{bgkfin}
\cite{apudir} The following classes have SAPU.
  \begin{enumerate}[(a)]   
 \item 
The class of posets with a maximum. 
\item
The class of directoids with a maximum.
\item
The class of bidirectoids with a maximum and a minimum.
\item
The class of  finite directoids.
\item
The class of  finite bidirectoids.
 \end{enumerate}
The following classes  have
SAP but not APU.
  \begin{enumerate}[(a)]
\setcounter{enumi}{5}
    \item    
 The classes of upward directed 
(downward directed, both upward and downward directed)
 posets.
\item
The class of directoids.
\item
The class of bidirectoids.
\end{enumerate} 
 \end{proposition}  

Using the methods of the present paper,
Proposition \ref{bgkfin} can be generalized so that 
structures with many directoid operations, with possibly many 
unary operations can be taken into account. 
We refer again to \cite{apudir} for more details. 

\subsubsection*{Order algebras.} \labbel{oasu} 
We now outline a general method to obtain 
classes of algebraic structures with SAPU from corresponding
classes of relational structures.
We first recall the original motivating example.

If $(P, \leq)$ is a poset, 
define a binary operation $\cdot$ 
(henceforth denoted by juxtaposition) on $P$ by
\begin{equation*}\labbel{oa} 
ab    =\begin{cases} a &   \text{if  $ a \leq b  $},\\ 
b  &    \text{otherwise.}
 \end{cases}
 \end{equation*}   
The structures $(P, \cdot)$
which can be obtained in this way 
have been described in \cite{Ng} under the name
pogroupoids, but   
are called 
\emph{order algebras} in the more recent
literature, e.~g.,  \cite{FJJMMM}. 
Clearly, we can retrieve the order $\leq$ from
$\cdot$ by setting $a \leq b$ if
$ab=a$.

Neggers \cite{Ng}
has showed, among other,  that
 if $P$, $Q$ 
are posets and $\iota: P \to Q$ is a function,
then $\iota$  is a homomorphism of the corresponding order algebras
if and only if $\iota$ is an ordermorphism with the further
property that
 any pair of incomparable elements are  sent either
to the same element, or to another pair of incomparable elements.
It follows immediately that if $\iota$ is injective, then
 $\iota$  is an embedding of order algebras
if and only if $\iota$ is an order-embedding,
since order-embeddings are exactly injective ordermorphisms
which send incomparable pairs to incomparable pairs.

\begin{corollary} \labbel{coroa}
The class of order algebras has SAPU.
\end{corollary} 

\begin{proof}
Immediate from the above observations and the fact, generalized here
in Proposition \ref{pu},  
that the class of posets has SAPU.
 \end{proof}  

\subsubsection*{Algebraizing relations.} \labbel{trf} 
The above argument has a more general flavor.
Let $\mathcal K$ be a class of structures such that 
$R(x,x,x, \dots)$ holds in every structure in $\mathcal K$,
for every relation symbol $R$
in the language. The language of $\mathcal K$
is allowed to contain 
relation, constant and function symbols.
Let us associate to $\mathcal K$ a class $\mathcal K_a$
defined in the following way. To every model $\mathbf A$ in $\mathcal K$
one associates a model $\mathbf A_a$  
 obtained from $\mathbf A$   
by replacing every $n$-ary  relation
$R$ 
by an $n$-ary  function $f_R$
defined by  
\begin{equation*}\labbel{agen} 
f_R(a_1,a_2,a_3, \dots)   =\begin{cases} a_1 &  
\text{if  $R(a_1,a_2,a_3, \dots) $ holds in $\mathbf A$},
\\[2pt] 
a_i  & 
 \begin{aligned}  &\text{otherwise, where $i$ is the smallest index}
\\[-4pt]
 &\text{such that $a_i \neq a_1$.}
\end{aligned}
 \end{cases}
 \end{equation*}   
Notice that in the second clause such an $i$ exists, 
since $R(a,a,a, \dots)$ holds in $\mathbf A$. 
The class $\mathcal K_a$ is the class of models
which can be obtained in this way.
 As in the case of order algebras,
from $\mathbf A_a$ we can retrieve the structure of $\mathbf A$.
It follows that the (class) function 
which sends $\mathbf A$ to $\mathbf A_a$
is bijective from $\mathcal K$ to $\mathcal K_a$.
 Again as in the case of order algebras, 
the notion of homomorphism in $\mathcal K_a$
is stronger than the notion of  homomorphism in $\mathcal K$;
however, the notions of embedding coincide. This is proved
just arguing
as in the last lines of \cite[proof of Theorem 1]{Ng}. 
Henceforth we get the following proposition.

\begin{proposition} \labbel{propa}
Under the above notations and conventions,
a class $\mathcal K$ has AP (APU, SAP, SAPU)
if and only if 
$\mathcal K_a$ has AP (APU, SAP, SAPU).

In particular, $\mathcal K_a$ has AP (and frequently SAPU)
for all the classes $\mathcal K$ described in
Propositions \ref{pu}, \ref{bgkfin}, \ref{dense}(a)-(b),
 Corollaries \ref{propj},  \ref{clp}(A)(B){\hspace{0 pt}}(C)(F)  and
Theorems \ref{due}(A)(B),  \ref{lopp}(1)-(6), \ref{mul}, \ref{gen}, \ref{cgk}.  
 \end{proposition} 

APU for order algebras and ``algebraized'' structures can be seen from a more
general perspective.

\begin{proposition} \labbel{fa}
Suppose that $\mathcal K$ is a class of structures for the same language.
  \begin{enumerate}    
\item 
Suppose further that,  
 for every $\mathbf A \in \mathcal K$,
every $n$-ary function $f _{ \mathbf A} $ 
and $a_1, \dots, a_n \in A$, we have
$f _{ \mathbf A}(a_1, \dots, a_n) \in \{ a_1, \dots, a_n \}  $.

If $\mathcal K$ has (S)AP and is closed
under taking substructures, then $\mathcal K$ has (S)APU.  
\item
More generally, 
suppose that,  
 for every $\mathbf A \in \mathcal K$,
every $n$-ary function $f _{ \mathbf A} $ 
and $a_1, \dots, a_n \in A$, there are a unary term
$t(x)$ and $i \leq n$  such that 
$f _{ \mathbf A}(a_1, \dots, a_n) = t(a_i)  $.

If $\mathcal K$ has (S)AP and is closed
under taking substructures, then $\mathcal K$ has (S)APU.  
\end{enumerate} 
 \end{proposition}  

\begin{proof} 
Any TBA triple $\mathbf A, \mathbf  B, \mathbf  C \in \mathcal K$
can be amalgamated to some $\mathbf  D \in \mathcal K$.
Under the assumptions, $A \cup B$ is a substructure of $\mathbf  D$.   
\end{proof}

\section{The joint embedding property into union.} \labbel{jep}

The classical joint embedding property (JEP)  \cite{H}
admits variations in the same spirit of Definition \ref{sapu}. 
For ordered sets with
operators, the corresponding theory is quite simple, in the sense that
the resulting disjoint embedding property into union
turns out to be generally either trivially true, or trivially false.
 In any case, JEP plays a fundamental role in model theory,
hence we shall explicitly mention the JEP-related properties.

\begin{definition} \labbel{jepu}
(a) Recall from Definition \ref{jepdef} that a   class $\mathcal K$ has the \emph{joint embedding  property}
(JEP) if, for every  
$\mathbf A, \mathbf B  \in \mathcal K$,
 there are a structure
$\mathbf D \in \mathcal K$ and  embeddings
$ \iota: \mathbf A \rightarrowtail \mathbf D$
and 
 $ \kappa: \mathbf B \rightarrowtail \mathbf D$. 
 
(b) $\mathcal K$ 
has the \emph{disjoint embedding  property}
(DJEP) if, under the assumptions from (a), $\mathbf D$, $\iota$
and $\kappa$  can be chosen in such a way that 
 the images of
$\iota$
and $\kappa$ are disjoint.

Clearly, this is impossible  if some constant 
is present in the language. 
This is the main difference with respect to SAP.

(c) We say that $\mathcal K$ 
has the \emph{joint embedding  property into union}
(JEPU) if, under the assumptions from (a), $\mathbf D$, $\iota$
and $\kappa$  can be chosen in such a way that 
 $D$ is the union of the images of
$\iota$
and $\kappa$.

(d) We say that $\mathcal K$ 
has the \emph{disjoint embedding  property into union}
(DJEPU)
if (b) and (c) can always be accomplished simultaneously, 
namely, $D$ can be chosen to be the
disjoint  union of the images of
$\iota$
and $\kappa$.
\end{definition}   

If $\mathcal K$ is closed under isomorphism,
then a remark parallel to \ref{iso} applies,
namely, $\mathcal K$ has DJEP if and only if, 
whenever  
$\mathbf A, \mathbf B  \in \mathcal K$
and $A \cap B = \emptyset $, then 
 there is a structure
$\mathbf D \in \mathcal K$
such that $\mathbf A \subseteq \mathbf  D$ 
and $\mathbf B \subseteq \mathbf  D$. Thus
$\mathcal K$ has DJEPU if and only if
$\mathbf  D$ as above can be chosen in such a way that
$D=A \cup B$.

The classical joint embedding property
dates back at least to \cite{F}. 
The disjoint embedding property
has been sometimes used in the literature, at least 
from the '80's in the last century, e.~g., Pouzet
\cite{P}.
For certain classes, e.~g.,
linearly or partially ordered sets,  
DJEPU is trivially satisfied
but, again, we do not know of a study of such properties
for their own sake.

If we  work in a language without constant
symbols, then in most arguments from the preceding sections
we may allow $\mathbf C$ to be empty.
As we mentioned, this special instance of AP (and  variations) turns out to be
exactly JEP (and  variations).
This is the reason why
remarks parallel to Observation \ref{obs}
and Propositions \ref{T}, \ref{un}, \ref{agg} hold,
possibly with a difference for languages with constant symbols.

\begin{observation} \labbel{obsj}
If some language $\mathscr L$ 
has no constant symbol, then the class $\mathcal K$ of all models for 
$\mathscr L$  has
DJEPU and any theory for $\mathscr L$ with only axioms of the form
\eqref{sent} from Proposition \ref{un}  has  DJEPU.

For arbitrary languages, Observation \ref{obs}(a)(b) and 
Propositions \ref{T}, \ref{un}(a){\hspace{0 pt}}(b), \ref{agg}(A1)-(A5) 
hold when 
(S)AP(U) is replaced everywhere by (D)JEP(U).
However, in Propositions \ref{T}(b) we need to assume that 
$\mathscr L'$ has no constant symbol. 
\end{observation}   

Notice that constants caused no trouble in getting
strong AP in Observation \ref{obs}(c), since constants
 are always already interpreted in 
$\mathbf C$.  On the other hand, as we have mentioned
in Definition \ref{jepu}(b),
the existence of some constant in the language always forbids
disjoint JEP (and frequently forbids JEP).

\begin{theorem} \labbel{frai}
All the classes and theories considered in 
Propositions \ref{pu}, \ref{bgkfin}(a)-(e),  
 Corollaries \ref{propj},   \ref{clp}(A)(B) and 
Theorems \ref{due}(A)(B), \ref{lopp}(1)-(4) , \ref{mul}, \ref{gen}, \ref{cgk} 
  have DJEPU,
with the provision that no new constant is added in 
the last two statements in Theorem \ref{gen}.

The classes considered in 
Corollary \ref{clp}(C)(F)
have JEPU
and the classes considered in 
Proposition \ref{bgkfin}(f)-(h) and 
Theorem \ref{lopp}(5)-(6) have
DJEP.

In particular, all the above classes have JEP.
 \end{theorem}  

\begin{proof}
 In most of the above arguments 
we have not forbidden the possibility that $\mathbf  C$ 
is empty, hence,
as remarked above, this special instance of
(S)AP(U)
provides (D)JEP(U). 
In fact, in this special case proofs turn out
to be generally much simpler: for example, 
in the case of posets, just consider 
$R=   R_{\mathbf A} \cup   R_{\mathbf B}$ 
on $A \cup B$. 
As another example, in the case of linearly ordered sets, 
just let all the elements of $\mathbf A$ to be $<$ than all the elements
of $\mathbf  B$. Adding operators presents no significant trouble, too. 

Some care is needed in the presence of constants.
However, in the case of bounded directoids, and variants,
the constants always generate isomorphic subalgebras,
hence we can consider, as $\mathbf  C$, any copy of this
``prime'' subalgebra, and then apply (S)AP(U).
\end{proof}

On the other hand, in comparison with 
Theorem \ref{apuvarthm},
 only varieties defined by a special kind of linear identities have
JEP, possibly, DJEPU. See
\cite{apuvar}, in particular,  
Theorem 3.1, Remark 4.1 and Corollary 4.6  there.

Dense linear orders provide interesting examples
concerning the  properties
dealt with in the present section.

\begin{proposition}   \labbel{dense}
(a) The theory of dense linear orderings has SAP,
DJEP, JEPU, APU 
but neither DJEPU nor SAPU. 

(a$'$) The theory of dense linear orderings
without endpoints has SAP,
DJEPU,  APU 
but  not SAPU.

(b)  The theory of dense linear orderings with a closure operation 
has SAP, DJEP, JEPU
but neither DJEPU, nor APU. 

(c)  The theory of dense linear orderings with two closure operations 
has not AP.
\end{proposition} 

 \begin{proof} 
(a) SAP follows immediately from the facts
that the class of linear orders has SAP (actually, SAPU)
and that every linear order can be embedded
into some dense linear order. 
DJEP is the special case when $\mathbf  C$ is the 
empty structure.

We now disprove SAPU. Let $\mathbf  C$ be $\mathbb Q$
with the standard order,
 consider two distinct copies 
$r$ and  $r'$ of the same real (not rational) number
and let $A= \mathbb Q \cup \{ r \} $,
 $B= \mathbb Q \cup \{ r' \} $.
However we linearly order $A \cup B$,
either $r'$ is the immediate predecessor 
of $r$, or conversely.
The counterexample works for
dense linear orders without endpoints, too. 

 Then we disprove DJEPU. 
Let $\mathbf A$ and $\mathbf  B$ be two disjoint copies of 
the real interval $[0,1]$ and suppose by contradiction
that $D=A \cup B$ can be densely linearly ordered   
extending the orders on $\mathbf A$ and $\mathbf  B$.
Define the following equivalence relation on $A$:
$ r \sim s$ if $\{ \, t' \in B \mid t' <_{ \mathbf  D} r  \,\}=
\{ \, t' \in B \mid t' <_{ \mathbf  D} s \,\}$.
Thus the $\sim$-equivalence classes partition $A$ 
and each   class is a convex subset of $A$, hence an interval.
It is easy to see that if $[0,1]$ is partitioned into intervals,
then at least one interval is a closed interval of the form $[r,s]$,
possibly with $r=s$.

Indeed, if the class of $0$ has the form $[0,v]$, we are done.
Otherwise, let $r$ be the largest real such that, for every 
$t < r$, the equivalence class of $t$ has the form $[u,v)$.
Namely, $r$ is the supremum  of those $v$
such that some equivalence class has the form 
 $[u,v)$ and also all preceding classes have that form.
 Then the class of $r$ has necessarily the form $[r,s]$,
since the form $[r,s)$ would contradict the definition of $r$. 

So let $[r,s]$ be a $ \sim$-equivalence class
and let 
$T'=\{ \, t' \in B \mid t' <_{ \mathbf  D} r  \,\}$.
If $T'$ is empty, then $0'$ in $B$ 
is the immediate successor of $s$  in $\mathbf  D$.
If $T'$ has a maximum $u'$ in $\mathbf  B$,
then $r$ is the immediate successor of $u'$
in $\mathbf  D$.
On the other hand, if the supremum $u'$ of $T'$ 
does not belong to $T'$, then
 $u'$ is the immediate successor of $s$
in $\mathbf  D$.
 In any case, we have found two elements
without intermediate elements,
 hence the order in $\mathbf  D$ is not dense, a contradiction.

Next, we prove JEPU. Actually, we show that DJEPU
fails ``for just one element'',
namely we get JEPU by identifying at most
one element from $\mathbf A$ with at most one element from $\mathbf  B$. 

So let $\mathbf A$ and $\mathbf  B$ be two dense linear
orderings. If either $\mathbf  A$ has no maximum
or $\mathbf B$ has no minimum, simply
put all the elements of $A$ before all the elements of $B$.
Otherwise, $\mathbf  A$ has a maximum $ \bar{a}$
and $\mathbf B$ has a minimum $ \bar{b}$.
Identify $ \bar{a}$ and $ \bar{b}$ in $\mathbf  D$ 
and, again,  put all the other
elements of $A$ before all the elements of $B$.

Having proved JEPU, it is rather easy to prove APU.

Let $\mathbf A_1$, $\mathbf  A_2$, $\mathbf  C$ 
be a TBA triple of dense linear orderings,
with embeddings $\iota_i : \mathbf  C \rightarrowtail \mathbf A_i$,
for $i=1,2$.  
Recall that if $\mathbf C$
is a linearly ordered set, a \emph{cut} of  $\mathbf C$
is a pair $(C', C'')$ such that 
$C' \cup C'' = C$
and $c' < c''$, for every     
$c' \in C'$ and $c'' \in  C''$.
We allow $C'$ or $C''$ to be empty.  
To any cut $(C', C'')$ of $\mathbf  C$ 
one associates on  $A_i \setminus \iota(C)$ the \emph{components}
$ \{ \, a \in A_i \setminus \iota(C) \mid \iota_i(c') <_ {\mathbf A_i}
 a <_ {\mathbf A_i} \iota_i(c''), \text{ for all } c' \in C'
\text{ and }   c'' \in C'' \,\}$, for $i=1,2$.
Conversely, to each
$a \in A_i \setminus \iota(C)$ one can associate the cut formed by  
$C' = \{ \,  c \in C \mid \iota_i(c) < _ {\mathbf A_i} a \, \} $
and 
$C'' = \{ \,  c \in C \mid a < _ {\mathbf A_i} \iota_i(c) \, \} $.
The nonempty components partition both 
$ A_1 \setminus \iota(C)$ and $A_2 \setminus \iota(C)$.
Moreover, if 
$(C', C'')$ is associated to $a \in A_i$
and $(C_*', C_*'')$ is associated to $a_* \in A_j$
with, say, $C' \subsetneq C_*'$,
then, for any possible amalgamating structure 
$\mathbf  D$ through embeddings
 $ \kappa _i : \mathbf  A_i \rightarrowtail \mathbf C$
($i=1,2$), we should have
$\kappa _i (a) <_ {\mathbf D} \kappa _j(a_*)$.  
This implies that, in order to construct 
$\mathbf  D$ and the $\kappa_i$'s, 
it is enough, for every cut,  to set the relative order 
between the elements of 
the components on 
$ A_1 \setminus \iota(C)$ and $ A_2 \setminus \iota(C)$ 
associated to the cut.  
See the proof of \cite[Thorem 3.1(a)]{lop} 
for more details.

We are almost done. 
If $E_i$, for $i=1,2$, are the components
 on $ A_i \setminus \iota(C)$
 associated to some cut, then
it is enough to embed the two $E_i$'s
into some dense linear order using JEPU.
It is easy to see that,  
letting the cut vary among all cuts of $\mathbf  C$
and putting together all the structures as above,
we get a dense linear order. 

(a$'$) As in (a), we have SAP since linear orders have SAP
and  every linear order can be embedded
into some dense linear order without endpoints.
The failure of SAPU has already been taken care of.
As far as APU is concerned, the argument in (a)
works in the present case, too, since if 
$\mathbf A_1$, $\mathbf  A_2$ and $\mathbf  C$ 
have no endpoint, then the model $\mathbf  D$ we have constructed 
has no endpoint. DJEPU is trivial, just let every element of $A$ 
be $\leq$ than every element of $\mathbf  B$.

(b) By Corollary \ref{clp}(B), the theory 
of linearly ordered sets with a closure operation $f$  has SAP.
Hence if $\mathbf A$, $\mathbf  B$, $\mathbf  C$ 
is a TBA triple of dense linear orders with a closure operation,
then there is an amalgamating linear order $\mathbf  E$
 with a closure operation. As an order, $\mathbf  E$
can be embedded into a complete dense linear order $\mathbf  D$
in such a way that, for every $d \in D$, there is
$e \in E$ such that  $d \leq _{\mathbf D} e$. 
Now define $f$ on $D$ by 
$f(d)= \inf  \{ \, f_{\mathbf E} (e)  \mid 
e \in E, \ d \leq _{\mathbf D} f_{\mathbf E} (e) \, \}  $
and it is easy to see 
that, with $f$ as defined,
$\mathbf  E$ embeds in $\mathbf  D$
and $f$ is a closure operation, thus $\mathbf  D$
amalgamates the original triple.  
We have proved SAP.

DJEP is the special case of
SAP when $\mathbf  C$ is the empty structure. 

We now check that JEPU holds.
As in (a), if either $\mathbf  A$ has no maximum
or $\mathbf B$ has no minimum, 
put all the elements of $A$ before all the elements of $B$.
Otherwise, $\mathbf  A$ has a maximum $ \bar{a}$
and $\mathbf B$ has a minimum $ \bar{b}$.
Then identify $ \bar{a}$ and  $ f(\bar{b})$
and put all the other elements of $A$ before all the
other elements of $B$. Thus all the (images of the) elements
of $A \setminus \{ \bar{a} \} $ precede
all the elements $b_1 \in B$
such that $b_1 < f(\bar{b})$ and all
such elements are bounded by 
$ \bar{a} = f(\bar{b})$. All the other elements of 
$B$ are larger. Notice that 
$f(\bar{a})\geq \bar{a}$ in $\mathbf A$, 
hence $f(\bar{a}) = \bar{a}$, since
$ \bar{a}$ is the maximum of $\mathbf A$.
Moreover, $f(f(\bar{b}))=f(\bar{b})$, 
since $f$ is a closure operation, hence the identification
 of $ \bar{a}$ and  $ f(\bar{b})$ is compatible.

Since  dense linear orderings have not DJEPU,
then  dense linear orderings with a closure operation have not DJEPU:
just consider the same counterexample with constant functions
added as operations.

In order to disprove APU, we shall modify  the counterexample
to SAPU given in (a). Fix $r \in \mathbb R \setminus \mathbb Q$
and  $q \in \mathbb Q$
with $r < q$.
Let $C= \mathbb Q$
with the standard order and  
define $f$ on   $C$ by
\begin{equation*}\labbel{cas}
 f(c)  =\begin{cases}
c  &    \text{if  $ c < r  $},
\\ 
q &    \text{if  $ r < c \leq q   $},
\\ 
c &    \text{if  $ q < c  $,}
\end{cases} 
\end{equation*} 
 thus $f$ is a closure operation on $C$. 

Let $A=C \cup \{ r \} $
with $f(r)=r$.
Let $r'$ be a copy of $r$
and define $\mathbf  B$ by 
$B=C \cup \{ r' \} $
with $f(r')=q$.
 As in (a),
if amalgamation into union holds, then
$r$ and  $r'$
should be identified, but this is impossible because of $f$.  

(c) Consider the  last example and add 
another operation $g$ defined as $f$
on $ \mathbb Q$ and such that 
$g(r)=q$ in $\mathbf A$ and 
$g(r')=r'$ in $\mathbf  B$. As above,
both  $f$ and  $g$ forbids the identification of 
$r$ and $r'$.  
Hence in any amalgamating structure
 we have either $r < r'$ or $r'<r$.
If  $r' < r$, then 
$q=f(r') \leq f(r)= r$, a contradiction. 
Symmetrically, $r'<r$ cannot hold,
thus AP fails.
\end{proof}

\section{More examples and counterexamples.} \labbel{coun} 

\begin{remark} \labbel{mor}
If in Definition \ref{sapu}
we replace embeddings with injective homomorphisms,
the results in the present note do not  necessarily hold.

(a) 
The class of posets has not  AP with respect to 
injective homomorphisms. Indeed, let 
$\mathbf C$ be the poset with just two incomparable elements
$ a$ and $b$. If $A= \{ a, b\} $ with $a < b$ in $\mathbf A$,
then the identity is an ordermorphism (but not an embedding!)
 $\iota _{\mathbf C, \mathbf A}$ from 
 $\mathbf C$ to $\mathbf A$.      
Similarly, let $B= \{ a, b\} $ with $b < a$ in $\mathbf B$
and $\iota _{\mathbf C, \mathbf B}$ be the identity map.
In any amalgamating poset  $\mathbf D$ we should have 
$ \iota _{\mathbf A, \mathbf D} (a) \leq \iota _{\mathbf A, \mathbf D} (b) $
and 
 $ \iota _{\mathbf B, \mathbf D} (b) \leq \iota _{\mathbf B, \mathbf D} (a) $.
Since 
we require
$\iota _{\mathbf C, \mathbf A} \circ \iota _{\mathbf A, \mathbf D}
=
\iota _{\mathbf C, \mathbf B} \circ \iota _{\mathbf B, \mathbf D}$,
then by antisymmetry
$ \iota _{\mathbf A, \mathbf D} (a) = \iota _{\mathbf A, \mathbf D} (b) $,
hence it is not possible to have
$ \iota _{\mathbf A, \mathbf D}$
an injective homomorphism.

(b) The main obstacle to  AP for injective homomorphisms 
in (a) is antisymmetry. 
In fact, the arguments in (a)
show that the class of sets with an 
antisymmetric binary relation has not AP
with respect to injective homomorphisms.

(c) In contrast with (a) and confirming (b), the  
class of preorders has  SAPU
with respect to  injective homomorphisms.

Parallel to Remark \ref{iso}, 
we can assume that $C=A \cap B$ 
and that the inclusions from 
$\mathbf C$ to $\mathbf A$ and 
from 
$\mathbf C$ to $\mathbf B$
are homomorphisms.
Then it is enough to endow  $A \cup B$
with the transitive closure of 
$\leq _{ \mathbf A} \cup \leq _{ \mathbf B}$. 
Of course, this is not the only possibility,
we could even  have done with
the \emph{discrete} preorder
(all pairs of elements are connected).

 (d) 
It is probably an interesting possibility to mix the two approaches, namely, to
consider embeddings (condition \eqref{R} is required)
with respect to  a certain set of relations, and homomorphisms 
 (condition \eqref{RR} is required)
with respect to  another set of relations.
Notice that there is no distinction between
embeddings and injective homomorphisms, when 
constants or functions are taken into account.

In the above proposal we intend to strictly remain
within the realm of model theory. The amalgamation property
can be defined  in a categorical setting and, of course,
this abstract setting encompasses all the above possibilities
\cite{KMPT}.
 \end{remark}

\begin{example} \labbel{1a}
(a) The theory (in the empty language)
 asserting that the universe has not cardinality $3$
has APU and SAP, but not SAPU.
Just let $|C|=1$ and $|A|=|B|=2$.

(b)
As above, the following theory $T$  in the empty language
has APU and SAP, but not SAPU.
The theory $T$ has sentences asserting
  \begin{enumerate}    
\item[($\diamondsuit$)] 
If there are at least $3$ distinct elements, then
there are at least $n$ distinct elements, 
  \end{enumerate}
for every $n\geq 3$.

The class of finite models of $T$
has APU but not SAP.

(c) If we consider the theory from (b)
in the language with a unary predicate $U$,
then $T$ has SAP, but 
the class of finite models of $T$
has not even AP.
Let $\mathbf  C$ have one element $c$, 
 let $\mathbf A$ have one more
element $a$ such that $U(a)$ and  
$\mathbf B$ have another
element $b$ such that not  $U(b)$.
Any amalgamating structure has at least 
$3$ elements, hence is infinite.  
 \end{example}

\begin{remark} \labbel{card}
(a) If we allow the empty model in the definition 
of SAPU and some class $\mathcal K$ has 
SAPU and both an empty model and 
a model of cardinality $1$,
then  $\mathcal K$ has models of any finite cardinality.

(b) If $\mathcal K$ has 
SAPU and has a model of cardinality $1$
which embeds into some model of cardinality $2$,
then  $\mathcal K$ has models of any finite nonzero cardinality.
 
(c) More generally,
if $\mathcal K$ has 
SAPU and has a model of cardinality $n$
which embeds into some model of cardinality $m > n$,
then  $\mathcal K$ has models of  cardinality
$m+k(m-n)$, for every $k$.
  \end{remark}

\begin{remark} \labbel{push}
(a) In most cases, the theories we have considered in this note
are universal Horn, hence they have pushouts.

In general, the homomorphisms
given by a pushout are not embeddings,
but if some class $\mathcal K$ has both AP and pushouts,
then, for every TBA triple in $\mathcal K$, the pushout 
is an amalgamating structure. In fact, in most cases, here
we have proved AP just by constructing the pushout
and showing that the homomorphisms towards the pushout
are embeddings.

(b) In the class $\mathcal K_0$ of sets
 (models for the empty language $\mathscr L_0$)
the pushout of a TBA triple $\mathbf A_0$, $\mathbf  B_0$
 and $\mathbf  C_0$ 
is the model $\mathbf  D_0$ over $D_0=A_0 \cup B_0$ without structure.
Hence some class $\mathcal K$ in some language $\mathscr L$ 
has SAP(U) if and only if, for every TBA triple
$\mathbf A$, $\mathbf  B$, $\mathbf  C$ in $\mathcal K$,
there is an amalgamating structure $\mathbf  D$ 
such that if $\mathbf  D^*$ is the pushout in $\mathcal K_0$
 of $\mathbf A _{ \restriction \mathscr L_0}  $,
 $\mathbf  B_{ \restriction \mathscr L_0} $ and
 $\mathbf  C_{ \restriction \mathscr L_0} $,
then $\mathbf  D^*$ is a subreduct (a reduct) of $\mathbf  D$.

(c) For transitive relations, a similar rephrasing
of the superamalgamation property is possible.

Let $\mathscr L_0= \{ R\} $ and $\mathscr L \supseteq \mathscr L_0$.  
If $\mathcal K$ is a class of models for $\mathscr L$
and $R$ is transitive in any member of $\mathcal K$,
then $\mathcal K$ has the superamalgamation property (into union) 
with respect to $R$ if and only if,
as above, 

(*) for every TBA triple
$\mathbf A$, $\mathbf  B$, $\mathbf  C$ in $\mathcal K$,
there is an amalgamating structure $\mathbf  D$ 
such that if $\mathbf  D^*$ is the pushout in $\mathcal K_0$
 of $\mathbf A _{ \restriction \mathscr L_0}  $,
 $\mathbf  B_{ \restriction \mathscr L_0} $ and
 $\mathbf  C_{ \restriction \mathscr L_0} $,
then $\mathbf  D^*$ is a subreduct (a reduct) of $\mathbf  D$.

(d) On the other hand, (*) form (c) above is not necessarily 
equivalent to the superamalgamation property, 
for relations which are  not supposed to be transitive.
Take $\mathscr L_0= \{ S \} $ and $\mathscr L= \{ R, S \} $.   
In the terminology of Theorem \ref{due}(A),
if $1 \in P$ and $1 \notin Q$, then $\mathcal K_{P,Q}$
has the superamalgamation property with respect to $S$,
but (*) is not satisfied, since, as shown by the proof,  we necessarily should add
$S$-related pairs  which are not related in the 
$\mathscr L_0$ pushout. 
 \end{remark}

We now show that, in order to get APU, we need to consider only unary
operations in  Proposition \ref{pu}(B)-(D).

\begin{proposition} \labbel{bin}
The theory $T$ of posets with a binary operation $f$  which is
order preserving on each component has not APU. 
 \end{proposition}    

\begin{proof}
Let $C= \{ c \} $ with the only possible structure.
Let $A= \{  a,c \} $ with $c < a$ and $f$ the projection onto the first
component.    
Let $B= \{  b,c \} $ with $c < b$ $b \neq a$ and $f$ the projection onto the second
component.    

In any amalgamating structure we should have
$f(a,b) \geq f(a,c)= a$ and 
$f(a,b) \geq f(c,b)= b$.
If  amalgamation is into union, 
then either $f(a,b) = a$, or 
$f(a,b) = b$.
Suppose the former,
hence necessarily $b \leq a$, since 
$b \leq f(a,b)$.  
Then
$c =f(c,a) \geq f(c,b) =b$, impossible.
\end{proof}

Not all possible variations on 
Proposition \ref{pu} hold. 

\begin{proposition} \labbel{rt}
\cite{apuduerel} Let $T$ be the theory  
with a binary reflexive and transitive relation $R$
and a unary function $f$  which strictly preserves $R$,
namely, 
\begin{equation*}    
\text{$d \mathrel { R } e $
and $d \neq e$ \quad imply \quad
both  
$f(d) \mathrel { R } f(e) $
and $ f(d) \neq f(e)$.}
 \end{equation*} 

Then $T$ has not AP. 
\end{proposition}

On the other hand, the theory $T$  
with a binary reflexive and transitive relation $R$
and a \emph{bijective} function $f$  which  preserves $R$
has superSAPU. This is immediate from the proof of 
Proposition \ref{pu} and also a special case of Theorem \ref{gen}.
Obviously, $R$-preserving bijective functions
are also strict  $R$-preserving.

The theory 
with an equivalence relation $R$
and a unary function $f$  which strictly preserves $R$
has SAP but not AP. See \cite{apuduerel} for details.

 \begin{example} \labbel{necess}
(a) We now show that the assumption that the $T_i$'s have SAPU 
in Proposition \ref{T}(a) cannot be weakened to APU,
even  for just  one among the $T_i$'s.

Let $T_1$ be the theory in the pure language of 
identity asserting that the universe has
cardinality $<3$.  Clearly, $T_1$
has APU. On the other hand,
SAP fails: just let $| C | = 1$
and $| A| = |B| = 2 $.

Let $T_2$ be the theory of 
partially ordered sets.
The classical proof that
$T_2$ has SAP actually provides SAPU,
as we noticed in Proposition \ref{pu}.  
 Let $C= \{ c \} $ 
 and $A= \{ c, a  \}$,  
$B= \{ b, c  \}$,
with 
$c < a$ in $\mathbf A$ and 
$b < c$ in $\mathbf B$.
The structures $\mathbf A$, $\mathbf B$ and  $\mathbf C$
are also models of $T_1$,
but any amalgamating structure 
must be of cardinality $\geq 3$,
hence is not a model of $T_1$.

(b) APU is not sufficient in Proposition \ref{T}(b), either.
Let $T_1$ be as above,  and $\mathscr L'= \{ U \} $,
where $U$ is a unary predicate. 
As above, let $| C | = 1$
and $| A| = |B| = 2 $.
Let $U(a)$ hold in $\mathbf A$, for   $ a \in A \setminus C$ 
and 
$U(b)$ fail in $\mathbf  B$, for   $ b \in B \setminus C$.
Thus in any amalgamating structure $\mathbf  D$ 
we have $a \neq b$, hence $| D| \geq 3$ and $\mathbf  D$ 
is not a model of $T_1$.   

In other words, $T_1$ has APU in the pure language of identity,
but has not even AP in the language  $\mathscr L'$. 
 \end{example}

A slightly more general version of 
Proposition \ref{T} holds with the same proof. 

\begin{proposition} \labbel{L}
 Suppose that $\mathscr L= \bigcup _ {i \in I}\mathscr L_i$
and the $\mathscr L_i$'s are pairwise disjoint  languages. 
Suppose that, for each $i \in I$,  $\mathcal K_i $
is a class of structures for $\mathscr L_i$ and
$\mathcal K_i $ has SAPU. 
Then $\mathcal K =  
\{ \mathbf A \mid \mathbf A \text{ is an $\mathscr L$-structure
and }
\mathbf A _{ \restriction \mathscr L_i}  \in \mathcal K_i,
\text{  for all } i \in I \} $ has SAPU.
 \end{proposition}

The next example is rather tricky, but it explains quite clearly
why the  ``S'' and the ``U'' in SAPU are necessary
in Proposition \ref{L}.

\begin{example} \labbel{***}
Take $\mathscr L_1 =\mathscr L_2 = \emptyset $ and
let $\mathcal K_1$ be the class of  models
of either  odd finite cardinality or of cardinality $\leq 2$.
Let $\mathcal K_2$ be the class of  models
of either  even finite cardinality or of cardinality $\leq 2$.
Both $\mathcal K_1$ and
$\mathcal K_2$ have SAP and APU but not SAPU. 
If
$\mathcal K$ is defined as in 
Proposition \ref{L}, then
$\mathcal K$ is the  class  of the models of cardinality $\leq 2$,
thus $\mathcal K$ has
 APU  but not SAP.

If $\mathscr L_3= \{ U \} $ and $\mathcal K_3$
is the class of all models for $\mathscr L_3$,
then, as in the proof in Example \ref{necess}(b), 
  $\mathcal K$ has not even AP.
 \end{example}

 \begin{example} \labbel{uneed}
We now provide counterexamples showing
that the version of Proposition \ref{un} fails
when (S)APU is weakened to (S)AP.

(a) Let $T$ be the theory of abelian groups
in the language with sum, opposite and a constant 
for the neutral element, with a further unary predicate $U$ 
and axioms stating, for every $n \in \mathbb N$:

\quad (i) if there are at least $n$ distinct elements such that $U(x)$,
then   there are at least $n$ distinct elements such that not  $U(x)$.

Clearly, (i) is expressible as a set of first-order
sentences.
The theory $T$ has SAP. Indeed,
given a TBA triple $\mathbf A$, $\mathbf  B$, $\mathbf  C$
of models of $T$, there is obviously
an amalgamating abelian group $\mathbf G$.
We have to interpret $U$ on $G$
in such a way that the expansion of 
$\mathbf G$ provides a model $\mathbf  D$ of $T$.
The interpretation of $U$ on $A \cup B$
is forced by the request that $\mathbf A \subseteq \mathbf  D$ 
and $\mathbf B \subseteq \mathbf  D$.
We can suppose that
 $\mathbf G$ properly extends both $\mathbf A$ and
$\mathbf  B$, since the other cases are trivial, Then,
considering laterals, we have $| G \setminus (A \cup B)| \geq |B \setminus A|$.
We are allowed to interpret $U$ in an arbitrary way
over $G \setminus (A \cup B)$, hence if we let $U(x)$
always fail on $G \setminus (A \cup B)$,
then $U(x)$ fails for at least half the elements  of $G \setminus A $.
Since $\mathbf A$ is a  model of $T$, then 
$U(x)$ fails for at least half the elements  of $A$.
In conclusion, with the above interpretation,
$U(x)$ fails for at least half the elements  of $G$,
thus $\mathbf  D$ is a model of $T$. 

Hence $T$ has SAP.
However, if 
$\sigma$ is the sentence 
$ \forall x \ (x=0 \vee U(x))$,
which has the form \eqref{sent},
then $T \cup \{  \sigma \} $
has not SAP. 
Indeed, let $\mathbf  C$ 
be a trivial group
in which $U(0)$ fails
and extend $\mathbf  C$ to  $\mathbf A$ and $\mathbf  B$, 
 two disjoint copies of $\mathbb Z_2$
in which $U(1)$ holds for both copies of $1$.
Any strong amalgamating group has cardinality
$\geq 4$, hence, if we interpret 
$U$  in such a way that (i) holds,
we have at least one element $d$ 
distinct from $0$ 
and such that not $U(d)$. 
But then $\sigma$ fails.

(b) The theory $T \cup \{  \sigma \} $ in the above counterexample
has not SAP, but $T \cup \{  \sigma \} $ has obviously AP.
Indeed, modulo isomorphism, the only models of $T \cup \{  \sigma \} $ are the trivial group and the two elements group, with $U$ interpreted as above.

However, the example can be modified in order to get 
 a theory $T$  with SAP such that $T \cup \{  \sigma \} $
has not even AP, for the same sentence $\sigma$ above. 

Simply consider the theory $T$ introduced in  (a), but in a language with a further unary
operation $f$ and no axiom mentioning $f$. The theory $T$ has SAP
even in the extended language: just
 amalgamate the structures without considering $f$
and then interpret $f$ in an arbitrary compatible way in the amalgamating structure
(this argument is the SAP-analogue of Proposition \ref{T}(b)).   
However, $T \cup \{  \sigma \} $, for $\sigma$ as in (a), has not AP: consider
the same counterexample as in (a), letting $f(1)=1$ in $\mathbf A$ and
$f(1)=0$ in $\mathbf  B$. By the considerations in (a),
the copies of $1$ in $\mathbf A$ and $\mathbf  B$ should be identified,
but this is prevented by the behavior of $f$.   

In fact, the above considerations are an example of a general phenomenon:
classes with SAP and classes with AP but not SAP are distinguished by their
behavior with respect  to expansions: the former classes are exactly
those classes with AP such that AP is preserved by expanding the language. 
We shall present details elsewhere.

(c) In the above examples the sentence
$\sigma$ is universal positive. We can modify the examples in such a way
that the sentence  is 
universal Horn. Let $T'$ be the theory of abelian groups
with a further unary predicate $V$ 
and axioms stating, for every $n \in \mathbb N$:

\quad 
(ii) if there are at least $n$ distinct elements such that not $V(x)$,
then   there are at least $n$ distinct elements such that   $V(x)$.

By the same arguments as in (a), $T'$ has SAP. 
Let $\sigma' $ be  
$ \forall x \ (V(x) \Rightarrow x=0 )$.
Then $T' \cup \sigma '$
has not SAP.
Let $\mathbf  C$ 
be a trivial group
in which $V(0)$ holds,
and let $\mathbf A$ and $\mathbf  B$ 
be two disjoint copies of $\mathbb Z_2$
in which $V(1)$ fails for both copies of $1$.
Then argue as above.

If we add a dummy unary function in the language, as in (b),
then $T'$
has still SAP in the expanded language,
while $T' \cup \sigma '$ has not AP. 
  
(d) The above example can be refined in order to obtain a
finitely axiomatizable 
universal theory $T''$ with SAP and a universal sentence $\sigma''$
such that $T'' \cup \sigma ''$ has not AP and $\sigma''$
has the following properties: 
$\sigma''$  is the universal closure of an atomic formula, 
  no constant appears in $\sigma''$ and
only one variable appears in $\sigma''$.

The language of  $T''$ consists
of a binary operation for addition, a unary operation for opposite,
two more unary operations $f$ and  $g$ and a unary predicate symbol
$V$. Axioms of $\mathcal T''$ contain axioms for abelian groups;
notice that we can do without a constant for the neutral element by asking
$(x+(-x))+y = y$, for all $x$ and  $y$.  
A further axiom of $T''$ asserts that 
$g$ is injective from $\{ \, x  \mid  \neg V(x) \,\}$   
 to $\{ \, x  \mid   V(x) \,\}$, namely,

(iii) for every $x,y$, if $x \neq y$, $\neg V(x)$
and $\neg V(y)$, then   
$g(x) \neq g(y)$, $ V(g(x))$
and $ V(g(y))$.

Since $g$ can be defined in an arbitrary
way on $\{ \, x  \mid   V(x) \,\}$, the arguments 
in (a)-(c) above show that
$T''$ has SAP. 
Let $\sigma''$ be $\forall x\ x+g(x)=x $.
We show that $T'' \cup \sigma ''$ has not AP.
As in (c), let 
 $\mathbf  C$ 
be a trivial group
in which (necessarily, because of (iii)) $V(0)$ holds.
Let $\mathbf A$ and $\mathbf  B$ 
be two copies of $\mathbb Z_2$
in which $V(0)$  holds and $V(1)$ fails, for both copies of $1$,
hence necessarily 
 $g(1)=0$. As in (b), let $f(1)=1$  
  in $\mathbf A$ and $f(1)=0$  
  in $\mathbf B$, thus the two copies of $1$
cannot be identified in any amalgamating structure. 
Hence APU fails and any amalgamating structure has
at least  one new element $d$,
hence at least one such  element.  
\end{example}

 \begin{example} \labbel{1}
In Proposition \ref{un}  it is necessary to assume that
in \eqref{sent} only one variable    
is bounded by $\forall$.
An example has been provided in Remark \ref{inf}(d); 
here is another example.

The theory $T$ without axioms in the language with two unary relations $R$ 
and $S$ has SAPU, by Observation \ref{obs}(c). 

If we add to $T$ the  axiom
\begin{equation}\labbel{rs}
\forall xy \ (R(x) \Rightarrow S(y)),
  \end{equation}    
then AP fails for the extended theory.

Take $C= \{ c \} $
and let $ \neg R(c)$ and
$S(c)$ hold in $\mathbf C$.

Let $\mathbf A$ over $A= \{ a,c \} $
extend $\mathbf C$ 
with $R(a)$ and $S(a)$.
Let $\mathbf B$ over $B= \{ b,c \} $
extend $\mathbf C$ 
with $\neg R(b)$ and $ \neg S(b)$.

Then $\mathbf A$ and $\mathbf B$
cannot be amalgamated over $\mathbf C$  
if we want that \eqref{rs} is satisfied. 
 \end{example}

 \begin{example} \labbel{injbij}
As we mentioned in Remark \ref{inf}(b),
if $\mathcal K$ is a class with SAPU in a language with a
unary  function
symbol $f$, then the subclass of $\mathcal K$ consisting of those
structures in which $f$ is bijective has SAPU.

We show that the corresponding statement is not true
when SAPU is weakened to  APU.
Let $T$ be the following theory in a language with 
two unary predicates $U$, $V$ and a unary function $f$.
The theory $T$ asserts that
  \begin{enumerate}    
\item
there is at most one element $x$ such that $V(x)$, and
\item 
for every $n \in \mathbb N$, 
if $V(x)$, then the elements
$x, f(x), \dots, f^n(x)$ are pairwise distinct
and do not lie in  $U$.   
  \end{enumerate}
 
Arguing as in Observation \ref{obs}(c)
we get that $T$ has APU, since, once there is some element $c$
such that $V(c)$, the theory describes completely 
the set $\{c, f(c), \dots, f^n(c), \dots\}$
and tells nothing about all the other potential elements,
except that $V$ never holds there.
In more detail, if $V(c)$, for some $c \in C$, 
then we have an amalgamating structure on $A \cup B$, as usual. 
This is the case also when $V(a)$, for some $a \in A$, but
not $V(b)$, for every $b \in B$, and conversely.
Similarly, we can amalgamate on the union when
not $V(x)$, for $x \in A \cup B$.  

On the other hand, SAP fails, since if in $\mathbf  C$
there is no element $x$ such that $V(x)$
but such elements exist both in $\mathbf A$ and $\mathbf  B$,
then they should be identified.
We still retain APU, since if $V(a)$ and  $V(b)$,
for some $a \in A$ and $b \in B$,
then we get an amalgamating structure by identifying     
$a$ with $b$,   $f(a)$ with $f(b)$,\dots, with no
further identification on $\mathbf  D$.  The identification can be
made coherently because of clause (2).

Let $T' = T \cup \{  \sigma \} $, where $\sigma$ 
says that $f$ is bijective.
Then   $T'$ has not AP.
Indeed, let 
$C= \{ c \} $ 
with $f(c)=c$, not $U(c)$
and not $V(c)$.   
Extend $\mathbf  C$ to models  $\mathbf A$ and $\mathbf  B$ 
by adding in each case a copy of $\mathbb Z  $,
with  $f(z)=z+1$, for $z \in \mathbb Z$, and  $V(1)$, in both cases,
but with $U(0)$  in $\mathbf A$ and not $U(0)$  in $\mathbf  B$.
By (1) the two copies of $1$   should be identified
in any amalgamating structure, hence, if $f$ is injective, then the two copies of $0$
should be identified, but this is impossible because of $U$. 
 \end{example}

The next example contrasts case (A3)
in Proposition \ref{agg}.

\begin{example} \labbel{exch}
Consider the  theory $T$ in the language with two binary relations
$R$ and $S$ and asserting that 
$ x \mathrel { R } y $ and   $ y \mathrel { R } z$
imply $ x \mathrel { S } z $.

(a) The theory $T$ has superSAPU with respect to $R$.

Indeed, given a TBA triple $\mathbf A$, $\mathbf  B$, $\mathbf  C$,
let   $R= R _{\mathbf A} \cup R _{\mathbf B}$
on $D= A \cup B$. Let $S$ be defined on $D$ 
by $d \mathrel { S  } e $ if 
 \begin{multline*}\labbel{alssu}   
\text{ \underline{either} $d,e \in A$
and   $d \mathrel {  S _{\mathbf A}} e $,
\underline{or} 
 $d,e \in B$
and   $d \mathrel {  S _{\mathbf B}} e $, \underline{or}}
\\
\text{$d \mathrel {  R } f \mathrel { R } e $,
for some $f \in D$.}
  \end{multline*} 
Then $\mathbf  D$ is a model of $T$ by construction.
It remains to show that $\mathbf  D$ extends $\mathbf A$ and $\mathbf  B$.
Again, this holds by construction, as far as $R$ is concerned;
moreover, the inclusions from $\mathbf A$ and $\mathbf  B$ to 
$\mathbf  D$ are homomorphisms.

Suppose that $d \mathrel { S } e$ in $\mathbf  D$ and, say,
$d,e \in A$. Then either $d \mathrel {  S _{\mathbf A}} e $,
or $d \mathrel {  R } f \mathrel { R } e $,
for some $f \in D$. By the definition of 
$R$ on $D$, we necessarily have $f \in A$
and  $d \mathrel {  R _{\mathbf A} } f \mathrel { R  _{\mathbf A}} e $,
thus $d \mathrel {  S _{\mathbf A}} e $, since $\mathbf A$ is a
model of $T$.
From $d,e \in A$ and $d \mathrel { S } e$ in $\mathbf  D$
we have got $d \mathrel {  S _{\mathbf A}} e $ in each case, and
this means that the inclusion of $A$ into $D$ is an embedding.
The argument for the inclusion of $B$ into $D$
 is symmetrical.

We have showed that $T$ has SAPU.
SuperSAPU with respect to $R$ follows from the definitions.

(b) There is a TBA triple $\mathbf A$, $\mathbf  B$, $\mathbf  C$ 
of models of $T$ such that no infinite chain of elements such that 
$a_1 \mathrel { S }a_2 \mathrel { S }a_3 \mathrel { S }a_4 \dots $ 
exists in $\mathbf A$, $\mathbf  B$ or $\mathbf  C$, but such a chain
exists in any amalgamating structure. 

Let $  C = \{ c_1, c_2, c_3, \dots  \} $ be a countably infinite set 
 with no pair of elements $R$-related and
no pair of elements $S$-related.
Let $A=C \cup \{ a_1, a_2, a_3, \dots  \}$ 
with $a_{i} \mathrel { R} c_{i} \mathrel { S} a_{i+1} $,
 for all indices $i$ 
and no other pair of elements $R$-
or $S$-related.
Let $B=C \cup \{ b_1, b_2, b_3, \dots  \}$ 
with $c_i \mathrel { R} b_i  \mathrel { S} c_{i}$, for all indices $i$,
and no other pair of elements $R$-
or $S$-related.
In any amalgamating structure
which is a model of $T$ we must have
$a_{1} \mathrel { S} b_{1} \mathrel { S} c_{1} 
\mathrel { S} a_{2} \mathrel { S} b_{2} \dots$

Thus, in contrast with the case
of partially ordered sets, 
Proposition \ref{agg}(A3),
  SAPU does not prevent the creation
of infinite chains of elements related
by the same relation.
\end{example}

\section{Further remarks.} \labbel{furt}

\begin{proposition} \labbel{senzau}
Suppose that $( T_i) _{i \in I} $
is a sequence of theories in disjoint languages $\mathscr L_i$. 
If each $T_i$ has SAP,
then $ T= \bigcup_{i \in I} T_i  $ has SAP. 
 \end{proposition}

\begin{proof}
(Sketch)
First observe that if $\mathcal K$ is 
a class of structures closed under
isomorphism with SAP and
 $\mathbf A \in \mathcal K$  
has no proper extension in $\mathcal K$,
then $\mathbf A$ has no proper substructure  
in $\mathcal K$.
Indeed, were $\mathbf  C$ a proper substructure of $\mathbf A$,
we could amalgamate $\mathbf A$ with an isomorphic copy
of $\mathbf A$ 
intersecting $A$ in $C$, 
getting a proper extension of $\mathbf A$.
Notice that here it is fundamental to assume the strong
version of AP.

Thus structures without  
proper extensions for some $T_i$  give no trouble.
Otherwise,
if  some model $\mathbf A$ of $T_i$ 
has a proper extension satisfying $T_i$,
then, for every infinite cardinal $\lambda \geq |\mathscr L_i|$,
  $\mathbf A$  has a proper extension 
of  cardinality $\lambda$ satisfying $T$.
This is immediate from the
 L\"owenheim-Skolem-Tarski Theorem 
if $\mathbf A$ is infinite; otherwise, by the preceding paragraph, 
 $T_i \cup Diag (\mathbf A)$
has  models of cardinality $\geq n$,
for arbitrary $n \in \mathbb N$, hence an infinite model
by compactness.

Thus, given a triple to
be amalgamated,
we can amalgamate their reducts
to the language of $T_i$
 to models 
$\mathbf  D_i$ such that  
$D_i \setminus (A \cup B)$
has  the same cardinality for each $i$. 
Since the languages are disjoint and 
each $T_i$ has SAP, 
we can arrange things in such a way that
the $\mathbf  D_i$'s have the same domain. 
\end{proof}
 
See \cite{BGR} for a detailed proof of a slightly
different statement, connections  with
quantifier-free interpolation and application to
verification and automated reasoning.

For theories with superAPU Propositions \ref{T}
and \ref{senzau}  
allow a generalization to non-disjoint languages.

\begin{proposition} \labbel{usup}
Suppose that $( T_i) _{i \in I} $
is a sequence of theories in  languages $\mathscr L_i$
and suppose that $\mathscr L_i \cap \mathscr L_j = \{ R \} $,
for $i \neq j \in I$, 
where $R$ is a binary relation symbol.

If each $T_i$  has superSAPU
and asserts that $R$ is transitive, then $T=\bigcup_{i \in I} T_i$
has superSAPU.
 \end{proposition} 

 \begin{proof}
Suppose that $\mathbf A$, $\mathbf  B$, $\mathbf  C$ 
is a TBA-triple of models of $T$ and, as usual, let
$D= A \cup B$.
Since $R$ is assumed to be transitive, if $\mathbf A$ and
$\mathbf  B$ embed in a model $\mathbf  D$ over $D$,
then the interpretation of $R$ is uniquely determined
by superSAPU. By assumption, each $\mathscr L_i$-reduct
of the triple can be superamalgamated  in some model over $D$.
 Then $R$ is interpreted
in the same way, for each $i \in I$. Since the languages pairwise
intersect in $ \{ R \} $, we can join all the interpretations
in a  model for the full language of $T$.  
 \end{proof} 

\begin{remark} \labbel{nectr}
The assumption that the common relation is transitive is necessary in 
Proposition \ref{usup}. Let $T_1$, resp., $T_2$  be the theories
of an antisymmetric relation $S$ with a finer partial order
$\leq_1$, resp, a finer partial order $\leq_2$.
Then both  $T_1$ and $T_2$
have superSAPU with respect to $S$, by Theorem \ref{due}(A).

On the other hand, $T_1 \cup T_2$ has not AP, by Proposition \ref{div+}(a). 
 \end{remark}

\begin{problem} \labbel{oth}
 Find other ways, besides Propositions \ref{T}, \ref{un},
\ref{agg}, \ref{L}, \ref{senzau} and   
Theorems \ref{mul} and \ref{gen}  
to merge classes with (S)AP
(not necessarily into union)
in such a way that a class
with (S)AP is obtained.
 \end{problem}

\begin{problem} \labbel{!!}
If $T$ is a first-order theory,
let $T_1$ be the set of all consequences of $T$
of the form \eqref{sent}.
By Proposition \ref{un}(c),
$T_1$ has SAPU, hence SAP.   
 
Due to the importance of SAP,
it is probably interesting to study the relationships between
$T$, $T_1$ and their models.

Under suitable assumptions, if the class $\mathcal K$ of finite substructures of 
models of $T_1$  
has JEP (this happens, for example, if the language of $T$ 
has no constant) then $\mathcal K$ has a Fra\"\i ss\'e limit
$\mathbf M$. See \cite[Section 7.1]{H}.  

Study the relationships between $T$ and $\mathbf M$.

Are there some other ways to extract a subtheory of $T$
having SAP and JEP?
\end{problem}

As another proposal for further research, 
it seems that SAPU fits well with 
models enriched with topological structures.
Moreover, SAPU can be used in order
to prove that certain theories have SAP,
though not necessarily SAPU.
An example appears  in the proof of 
\cite[Theorem 3.5]{J},
where SAPU for posets is 
implicitly used in order to 
prove SAP for lattices. 
Hence it is likely that  the present methods can be extended 
in order to prove SAP for many more theories.

\end{document}